\documentclass[5p]{elsarticle}

\usepackage{lineno,hyperref}
\modulolinenumbers[5]

\journal{Systems \& Control Letters}

\usepackage{amssymb}
\usepackage{amsmath}

\usepackage{graphicx}
\begin{document}

\begin{frontmatter}
\title{Constructing Artificial Traffic Fluids by Designing Cruise Controllers\tnoteref{mytitlenote}}
\tnotetext[mytitlenote]{The research leading to these results has received funding from the European Research Council under the European Union’s Horizon 2020 Research and Innovation programme/ ERC Grant Agreement n. [833915], project TrafficFluid.}

\author{Iasson Karafyllis }
\address{Dept. of Mathematics, National Technical University of Athens, Zografou Campus, 15780, Athens, Greece.}
 \author{Dionysis Theodosis${^1}$} \author{Markos Papageorgiou${^{1,2}}$}
\address{1. Dynamic Systems and Simulation Laboratory,  Technical University of Crete, Chania, 73100, Greece\\2. Faculty of Maritime and Transportation, Ningbo University, Ningbo, China.\vspace*{-3em}}

\begin{abstract}
In this paper, we apply a Control Lyapunov Function methodology to design two families of cruise controllers for the two-dimensional movement of autonomous vehicles on lane-free roads using the bicycle kinematic model. The control Lyapunov functions are based on measures of the energy of the system with the kinetic energy expressed in ways similar to Newtonian or relativistic mechanics. The derived feedback laws (cruise controllers) are decentralized, as each vehicle determines its control input based on its own speed and on the relative speeds and distances from adjacent vehicles and from the boundary of the road. Moreover, the corresponding macroscopic models are derived, obtaining fluid-like models that consist of a conservation equation and a momentum equation with pressure and viscous terms. Finally, we show that, by selecting appropriately the parameters of the feedback laws, we can determine the physical properties of the ``traffic fluid'', i.e. we get free hand to create an artificial fluid that approximates the emerging traffic flow.
\end{abstract}

\begin{keyword}
Control Lyapunov Function, Cruise Controllers, Artificial Traffic Fluid, Lane-Free Roads
\end{keyword}

\end{frontmatter}


\section{Introduction}
Traffic flow theory includes two fundamental classes of mathematical models: microscopic and macroscopic. Microscopic traffic models aim to describe the longitudinal (car-following) and lateral (e.g. lane-changing) movement of individual vehicles in the traffic stream; while macroscopic traffic flow models reflect the collective behavior of vehicles by use of aggregate variables (flow, density, and mean speed of vehicles).

The mathematical description of conventional traffic through microscopic models has been studied extensively with various contributions and applications. Recent advances in technology have revolutionized vehicle automation with different kinds of driver support systems (see for instance \cite{9}, \cite{18}, \cite{22}). In the era of connected and automated vehicles, new perspectives and principles have also been suggested \cite{19}, where autonomous vehicles can move on the two-dimensional surface of lane-free roads (\cite{11}, \cite{31}) without abiding to lane discipline, something that may improve traffic flow and increase capacity of highways and arterials. 

Macroscopic traffic flow modelling started in the 1950s and continued to this day with a variety of models and approaches, see for instance \cite{Bek}, \cite{3}, \cite{5}, \cite{7}, \cite{10}, \cite{12}, \cite{32}, \cite{33} and references therein. Furthermore, several methodologies have been suggested to derive macroscopic models from microscopic models, see for instance \cite{3}, \cite{4}, \cite{7}, \cite{8}, \cite{29}  and references therein.

In this paper, we extend the Control Lyapunov Function (CLF) methodology presented in \cite{11} to derive families of cruise controllers for autonomous vehicles on lane-free roads (Section 3). The CLFs also act as size functions (Lemma 1 and Lemma 2) guaranteeing that the closed-loop system is well-posed. The construction of the CLF is based on measures of the total energy of the system. By expressing the kinetic energy in ways similar to Newtonian or relativistic mechanics, two respective families of cruise controllers are obtained that satisfy the following properties globally (Theorem 1 and Theorem 2): (i) there are no collisions among vehicles or with the boundary of the road; (ii) the speeds of all vehicles are always positive and remain below a given speed limit; (iii) the speeds of all vehicles converge to a given speed set-point; and (iv) the accelerations, lateral speeds, and orientations of all vehicles tend to zero.  The proofs of the above results can be found in Section 5. The proposed families of cruise controllers are decentralized (per vehicle) and require either the measurement only of the distances from adjacent vehicles (inviscid cruise controllers) or the measurement of speeds of and distances from adjacent vehicles (viscous cruise controllers). In contrast to the analysis presented in \cite{11}, the vehicles are not assumed to have same size, i.e. each vehicle may have different size.

Finally, using the methodology  in \cite{7}, we formally derive the macroscopic models that correspond to the closed-loop systems with the derived cruise controllers (Section 4). The resulting macroscopic models are very similar to models describing the flow of a Newtonian, compressible fluid in a porous or non-porous medium. We also provide the explicit formulae that relate the physical characteristics of the ``traffic fluid'' to the parameters of the cruise controllers. This implies that, by changing the functions and the parameters of the cruise controllers, we can actually determine the physical characteristics of the ``traffic fluid'', i.e. \textit{{we get largely free hand to create an artificial fluid}} that approximates the emerging traffic flow. 

To understand how far the implications of the relations between the cruise controller parameters and the characteristics of the traffic fluid go, it is important to notice that, for isentropic (or barotropic) flow of gases, the dynamic viscosity and the pressure are always increasing functions of the fluid density (see the discussion in \cite{13}). However, using the proposed families of cruise controllers, it is possible to obtain a traffic fluid with dynamic viscosity and pressure that are non-monotone functions of the fluid density (and can have local minima). Thus, the traffic fluid can have very different physical properties from the properties of real compressible fluids (mainly gases). \textit{{Therefore, it can be claimed that the cruise control design procedure becomes the design procedure (with many degrees of freedom) of an artificial fluid}} This is a mathematically founded realisation of the incentive, expressed in \cite{19}, to design future traffic flow as an efficient artificial fluid via appropriate design of the underlying vehicle movement strategies.

\textbf{Notation.} Throughout this paper, we adopt the following notation. By $\mathbb{R} _{+} :=[0,+\infty )$ we denote the set of non-negative real numbers. By $|x|$ we denote both the Euclidean norm of a vector $x\in \mathbb{R} ^{n} $ and the absolute value of a scalar $x\in \mathbb{R} $. By $x'$ we denote the transpose of a vector $x\in \mathbb{R} ^{n} $. By $dist(x,A)=\inf \left\{|x-y|:y\in A\right\}$ we denote the Euclidean distance of the point $x\in \mathbb{R} ^{n} $ from the set $A\subset \mathbb{R} ^{n} $. Let $A\subseteq \mathbb{R} ^{n} $ be an open set. By  $C^{0} (A,\Omega )$, we denote the class of continuous functions on $A\subseteq \mathbb{R} ^{n} $, which take values in $\Omega \subseteq \mathbb{R} ^{m} $. By $C^{k} (A;\Omega )$, where $k\ge 1$ is an integer, we denote the class of functions on $A\subseteq \mathbb{R} ^{n} $  with continuous derivatives of order $k$, which take values in $\Omega \subseteq \mathbb{R} ^{m} $. When $\Omega =\mathbb{R} $ the we write $C^{0} (A)$ or $C^{k} (A)$.  For a function $V\in C^{1} (A\, ;\, \mathbb{R} )$, the gradient of $V$ at $x\in A\subseteq \mathbb{R} ^{n} $, denoted by $\nabla V(x)$, is the row vector $\left[\frac{\partial V}{\partial x_{1} } (x)\, \, \cdots \, \, \frac{\partial V}{\partial x_{n} } (x)\right]$. Let $a(n),b(n)$ be quantities depending on the integer $n\ge 1$. We say that $a(n)=b(n)+O(n^{-p} )$, where $p\ge 1$, if there exists a constant $K>0$ (independent of $n\ge 1$) such that $\left|a(n)-b(n)\right|\le K\, n^{-p} $ for all $n\ge 1$.

\section{Vehicular Model Description }
 
Consider $n$ vehicles on a lane-free road of constant width $2a>0$, where the movement of each vehicle $i\in \{ 1,...,n\} $ is described by the model  
\begin{equation} \label{GrindEQ__2_1_} 
\dot{x}_{i} =v_{i} \cos (\theta _{i} ),\,
\dot{y}_{i} =v_{i} \sin (\theta _{i} ) \,\,
\dot{\theta }_{i} =\sigma _{i}^{-1} v_{i} \tan (\delta _{i} ) ,\,
\dot{v}_{i} =F_{i}
\end{equation} 
where $\sigma _{i} >0$ is the length of vehicle $i$ (a constant). Here, $(x_{i} ,y_{i} )\in \mathbb{R} \times \left(-a,a\right)$ is the reference point of the $i$-th vehicle in an inertial frame with Cartesian coordinates $(X,Y)$, with $i\in \{ 1,...,n\} $ and is placed at the midpoint of the rear axle of the vehicle; $v_{i} \in (0,v_{\max } )$ is the speed of the $i$-th vehicle at the point $(x_{i} ,y_{i} )$, where $v_{\max } >0$ denotes the road speed limit; $\theta _{i} \in \left(-\frac{\pi }{2} ,\frac{\pi }{2} \right)$ is the heading angle (orientation) of the $i$-th vehicle with respect to the $X$ axis; $\delta _{i} $ is the steering angle of the front wheels relative to the orientation $\theta _{i} $ of the $i$-th vehicle; and $F_{i} $ is the acceleration of the $i$-th vehicle. Model \eqref{GrindEQ__2_1_} is known as the bicycle kinematic model (see \cite{22}), and has been used to represent vehicles due to its simplicity to capture vehicle motion. In what follows, we use the input transformation $\delta _{i} :=\arctan \left(\sigma _{i} v_{i}^{-1} u_{i} \right)$, to obtain the following system 
\begin{equation} \label{GrindEQ__2_2_} 
\begin{aligned}
\dot{x}_{i} =&v_{i} \cos (\theta _{i} ) ,\;\,
\dot{y}_{i} =v_{i} \sin (\theta _{i} ),\;\,
\dot{\theta }_{i} =u_i ,\;\,
\dot{v}_{i} =F_{i} \end{aligned} 
\end{equation} 
for $i=1,...,n$, where $u_{i} $, $F_{i} $, are the inputs of the system.  \nolinebreak

Let $v^{*} \in (0,v_{\max } )$ be a given speed set-point and define the set 
\begin{equation} \label{GrindEQ__2_3_} 
S:=\mathbb{R} ^{n} \times \left(-a,a\right)^{n} \times \left(-\varphi ,\varphi \right)^{n} \times (0,v_{\max } )^{n}  
\end{equation} 
where $\varphi \in \left(0,\frac{\pi }{2} \right)$ is an angle that satisfies
\begin{equation} \label{GrindEQ__2_4_} 
\cos \left(\varphi \right)>\frac{v^{*} }{v_{\max } }.                  
\end{equation} 
The set $S$ in \eqref{GrindEQ__2_3_} describes all possible states of the system of $n$ vehicles. Specifically, each vehicle should stay within the road, i.e., $(x_{i} ,y_{i} )\in \mathbb{R} \times \left(-a,a\right)$ for $i=1,...,n$; moreover, the vehicles should not be able to turn perpendicular to the road, hence it should hold that $\theta _{i} \in \left(-\varphi ,\varphi \right)$ for $i=1,...,n$. The constant $\varphi $ can be understood as a safety constraint, which  restricts  the movement of a vehicle; finally, the speeds of all vehicles should always be positive, i.e., no vehicle moves backwards; and all vehicles should respect the road speed limit. 

We define the distance between vehicles by
\begin{equation}\label{GrindEQ__2_5_}\hspace{-0.08em}
d_{i,j} :=\sqrt{(x_{i} -x_{j} )^{2} +p_{i,j} (y_{i} -y_{j} )^{2} },\textrm{ for }  i,j=1,...,n
\end{equation}
where $p_{i,j} \ge 1$ are weighting factors that satisfy $p_{i,j} =p_{j,i} $ for all $i,j=1,...,n$. Note that, for $p_{i,j} =1$, we obtain the standard Euclidean distance, while for larger values of $p_{i,j} >1$, we have an ``elliptical'' metric, which can approximate more accurately the dimensions of a vehicle. For the case of $n$ vehicles of equal length, the optimal selection of a single $p$ can be found in \cite{11}.
Let
\begin{equation} \label{GrindEQ__2_6_} 
w=(x_{1} ,...,x_{n} ,y_{1} ,...,y_{n} ,\theta _{1} ,...,\theta _{n} ,v_{1} ,...,v_{n} )'\in \mathbb{R} ^{4n} .                                         
\end{equation} 
Due to the various constraints explained above, the state space of the model \eqref{GrindEQ__2_2_} is the set
\begin{equation} \label{GrindEQ__2_7_} 
\Omega :=\left\{\, w\in S\, \, :\, \; d_{i,j} >L_{i,j} ,\; i,j=1,...,n\, ,\, j\ne i\, \right\} 
\end{equation} 
where $L_{i,j} $, $i,j=1,...,n$, $i\ne j$, are positive constants and represent the minimum distance between a vehicle $i$ and a vehicle $j$, with $L_{i,j} =L_{j,i} $ for $i,j=1,...,n$, $i\ne j$. To have a well-posed closed-loop system on the state space $\Omega \subset \mathbb{R} ^{4n} $, the control inputs $u_{i} $ and $F_{i} $, $i=1,...,n$, should be given by appropriate feedback laws which are designed in such a way that every solution of \eqref{GrindEQ__2_2_} satisfies the following implication: $w(0)\in \Omega \Rightarrow w(t)\in \Omega $ for all $t\ge 0$ (see Section 3).

Finally, it should be noted that model \eqref{GrindEQ__2_2_} is nonlinear not only because of the nonlinearities appearing on the right-hand sides of \eqref{GrindEQ__2_2_}, but also due to the fact that the state space $\Omega $ is not a linear subspace of $\mathbb{R} ^{4n} $, but an open set (see  \cite{25} for the extension of the Input-to-State Stability property to systems defined on open sets). As noted in Section 1, model \eqref{GrindEQ__2_2_} with state space given by \eqref{GrindEQ__2_7_} is an extension of the model given in \cite{11}, where all vehicles were assumed to be identical and all distances between vehicles were given by the same value of $p$ in \eqref{GrindEQ__2_5_}.

\section{Two Families of Cruise Controllers}

\subsection{Preliminaries}

In this section we present two families of cruise controllers for vehicles operating on lane-free roads that satisfy the following properties:
\begin{itemize}
\item[(P1)]  Well-posedness requirement: For each $w(0)\in \Omega $, there exists a unique solution $w(t)\in \Omega $ defined for all $t\ge 0$. According to \eqref{GrindEQ__2_7_}, this requirement implies that there are no collisions between vehicles (since $d_{i,j} (t)>L_{i,j} $ for $t\ge 0$, $i,j=1,...,n$, $j\ne i$) or with the boundary of the road (since $y_{i} (t)\in (-a,a)$ for $t\ge 0$); the speeds of all vehicles are always positive and remain below the given speed limit (since $v_{i} (t)\in (0,v_{\max } )$ for all $t\ge 0$); and the orientation of each vehicle is always bounded by the given value $\varphi \in \left(0,\frac{\pi }{2} \right)$ (since $\theta _{i} (t)\in (-\varphi ,\varphi )$ for $t\ge 0$).

\item[(P2)]  Asymptotic requirement: The orientation of each vehicle satisfies $\mathop{\lim }\limits_{t\to +\infty } \left(\theta _{i} (t)\right)=0$ for $i=1,...,n$, and the speeds of all vehicles satisfy $\mathop{\lim }\limits_{t\to +\infty } \left(v_{i} (t)\right)=v^{*} $, $i=1,...,n$, for a given a longitudinal speed set-point $v^{*} \in (0,v_{\max } )$. Moreover, the accelerations, angular speeds, and lateral speeds of all vehicles tend to zero, i.e., $\mathop{\lim }\limits_{t\to +\infty } \left(F_{i} (t)\right)=0$,  $\mathop{\lim }\limits_{t\to +\infty } \left(u_{i} (t)\right)=0$, and $\mathop{\lim }\limits_{t\to +\infty } \left(\dot{y}_{i} (t)\right)=0$, for $i=1,...,n$. 
\end{itemize}

Let $V_{i,j} :(L_{i,j} ,+\infty )\to \mathbb{R} _{+} $, $U_{i} :(-a,a)\to \mathbb{R} _{+} $, $i,j=1,...,n$, $j\ne i$ be $C^{2} $ functions and let $\kappa _{i,j} :(L_{i,j} ,+\infty )\to \mathbb{R} _{+} $,$i,j=1,...,n$, $j\ne i$  be $C^{1} $ functions that satisfy the following properties
\begin{align} 
\mathop{\lim }\limits_{d\to L_{i,j}^{+} }& \left(V_{i,j} (d)\right)=+\infty  \label{GrindEQ__3_1_} \\
V_{i,j} (d)&=0, \textrm{ for all } d\ge \lambda \label{GrindEQ__3_2_}\\
V_{i,j} (d)&\equiv V_{j,i} (d),\,\, i,j=1,...,n,j\ne i  \label{GrindEQ__3_3_} \end{align}\begin{align}
\mathop{\lim }\limits_{y\to (-a)^{+} }\left(U_{i} (y)\right)=& +\infty , \mathop{\lim }\limits_{y\to a^{-} } \left(U_{i} (y)\right)=+\infty  \label{GrindEQ__3_4_} \\
U_{i} (0)=&0 \label{GrindEQ__3_5_} \\
\kappa _{i,j} (d)=&0, \textrm{ for all } d\ge \lambda \label{GrindEQ__3_6_}\\
\kappa _{i,j} (d)\equiv &\kappa _{j,i} (d),\,\, i,j=1,...,n,j\ne i \label{GrindEQ__3_7_} 
\end{align} 
where $\lambda $ is a positive constant that satisfies
\begin{equation} \label{GrindEQ__3_8_} 
\lambda >\max \left\{L_{i,j} ,i,j=1,...,n,i\ne j\right\} .  
\end{equation} 
The families of functions $V_{i,j} $ and $U_{\iota } $ in \eqref{GrindEQ__3_1_}, \eqref{GrindEQ__3_2_}, \eqref{GrindEQ__3_3_}, and \eqref{GrindEQ__3_4_}, \eqref{GrindEQ__3_5_}, respectively, are potential functions, which have been used to avoid collisions between vehicles and road boundary violation (see for instance \cite{30}). Condition \eqref{GrindEQ__3_3_} implies that if a vehicle $i$ exerts a force to vehicle $j$, then vehicle $j$ exerts the opposite force to vehicle $i$. The functions $\kappa _{i,j} $ are used in the subsequent analysis for the introduction of a viscous-like behavior of the vehicles. 

We exploit next a Control Lyapunov Function (CLF) methodology and the potential functions above to derive families of cruise controllers for autonomous vehicles on lane-free roads that satisfy properties (P1) and (P2). The construction of the Lyapunov function is based on measures of the total energy of the system. Depending on how the kinetic energy of the system is expressed, we obtain two different families of cruise controllers. If the kinetic energy is expressed in a fashion similar to that of Newtonian mechanics, we call the corresponding controller a \textit{Newtonian Cruise Controller} (NCC); while, when the kinetic energy is expressed in terms similar to those of relativistic mechanics, we call the corresponding controller a \textit{Pseudo-Relativistic Cruise Controller} (PRCC). The main difference of those two approaches is that, in relativistic mechanics, the kinetic energy increases to infinity when an object's speed approaches the speed of light; while the kinetic energy in Newtonian mechanics continues to increase without bound as the speed of an object increases. 
Finally, when $\kappa _{i,j} (d)\equiv 0$ for $i,j=1,...,n,j\ne i$, we call the controller \textit{``inviscid''} since the corresponding macroscopic model does not contain a viscosity term; otherwise, the corresponding controller is called \textit{``viscous''}.

\subsection{Newtonian Cruise Controller (NCC)}

The CLF in this case is given by the formula
\begin{equation} \label{GrindEQ__3_9_} 
\begin{aligned} {H(w):} & {=}  {\frac{1}{2} \sum _{i=1}^{n}\left(v_{i} \cos (\theta _{i} )-v^{*} \right)^{2}  +\frac{b}{2} \sum _{i=1}^{n}v_{i}^{2} \sin ^{2} (\theta _{i} ) } \\ 
&+\sum_{i=1}^nU_i(y_i)+\frac{1}{2}\sum_{i=1}^n\sum_{i\neq j}V_{i,j}(d_{i,j})\\
&+A\sum_{i=1}^n\left(\frac{1}{cos(\theta_i)-\cos(\varphi)}-\frac{1}{1-\cos{\varphi}}\right)
\end{aligned}  
\end{equation} 
where $A>0$, $b>1-\frac{v^{*} }{v_{\max } }>0 $ (recall that $v^{*} \in (0,v_{\max } )$) are parameters of the controller and the Lyapunov function. The function $H$ in \eqref{GrindEQ__3_9_} is based on the total mechanical energy of the system of $n$ vehicles. Specifically, the first two terms ($\frac{1}{2} \sum _{i=1}^{n}(v_{i} \cos (\theta _{i} )-v^{*} )^{2}  +\frac{b}{2} \sum _{i=1}^{n}v_{i}^{2} \sin ^{2} (\theta _{i} ) $) are related to the kinetic energy of the system of $n$ vehicles  relative to an observer moving along the $x-$direction with speed equal to $v^{*} $(as in classical mechanics); they penalize the deviation of the longitudinal and lateral speeds from their desired values $v^{*} $ and zero, respectively. The sum of the third and fourth term ($\sum _{i=1}^{n}U_{i} (y_{i} )  +\frac{1}{2} \sum _{i=1}^{n}\sum _{j\ne i}$ $V_{i,j} (d_{i,j} )  $), which are based on the potential functions \eqref{GrindEQ__3_1_} and \eqref{GrindEQ__3_4_}, is related to the potential energy of the system. Finally, the last term of \eqref{GrindEQ__3_9_} ($A\sum _{i=1}^{n}(\frac{1}{\cos (\theta _{i} )-\cos (\varphi )} -\frac{1}{1-\cos (\varphi )} ) $) is a penalty term that blows up when $\theta _{i} \to \pm \varphi $.

While the CLF \eqref{GrindEQ__3_9_} has characteristics of a size function (see \cite{25}), it is not a (global) size function, since $H$ takes finite values for $v_{i} \notin (0,v_{\max } )$. The following lemma shows the partial size function properties of the Lyapunov function $H$.  

\textbf{Lemma 1: }\textit{Let constants $A>0$, $v_{\max } >0$, $v^{*} \in (0,v_{\max } )$, $L_{i,j} >0$, $i,j=1,...,n$, $i\ne j$, $\lambda >0$ that satisfies \eqref{GrindEQ__3_8_}, $\varphi \in \left(0,\frac{\pi }{2} \right)$ that satisfies \eqref{GrindEQ__2_4_}, and define the function $H:\Omega \to \mathbb{R} _{+} $ by means of \eqref{GrindEQ__3_9_}, where $\Omega $ is given by \eqref{GrindEQ__2_7_}. Then, there exist non-decreasing functions $\omega :\mathbb{R} _{+} \to \left[0,\varphi \right)$, $\eta _{i} :\mathbb{R} _{+} \to [0,a)$, $i=1,...,n$, and for each pair $i,j\in \left\{1,...,n\right\}$, $i\ne j$, there exist non-increasing functions $\rho _{i,j} :\mathbb{R} _{+} \to (L_{i,j} ,\lambda ]$ with $\rho _{i,j} (s)\equiv \rho _{j,i} (s)$, such that the following implications holds:}
\begin{equation} \label{GrindEQ__3_10_}
\begin{aligned}
w\in& \Omega \Rightarrow  \left|\theta _{i} \right|\le \omega \left(H(w)\right),\left|y_{i} \right|\le \eta _{i} \left(H(w)\right),\\& d_{i,j} \ge \rho _{i,j} \left(H(w)\right),\textrm{ for }i,j=1,...,n\, ,\, j\ne i.       
\end{aligned}
\end{equation} 

\noindent Based on the CLF \eqref{GrindEQ__3_9_}, we obtain the feedback laws  
\begin{align} 
&\begin{aligned}&{u_{i} }  {=}\\&  {\left(Z_{i} (w)-U'_{i} (y_{i} )-\sum _{j\ne i}p_{i,j} V'_{i,j} (d_{i,j} )\frac{(y_{i} -y_{j} )}{d_{i,j} }  -b\sin (\theta _{i} )F_{i} \right)} \\  &{\times \left(v^{*} +\frac{A}{v_{i} \left(\cos (\theta _{i} )-\cos (\varphi )\right)^{2} } +v_{i} \cos (\theta _{i} )(b-1)\right)^{-1} } \end{aligned} 
 \label{GrindEQ__3_11_}\\
&F_{i} =-\frac{1}{\cos (\theta _{i} )} \left(k_{i} (w)\left(v_{i} \cos (\theta _{i} )-v^{*} \right)+\Lambda _{i} (w)\right) \label{GrindEQ__3_12_}
\end{align} 
where 
\begin{align} 
&\begin{aligned}
Z_{i} (w):=&-\mu _{1} v_{i} \sin (\theta _{i} )\\&+\sum _{j\ne i}\kappa _{i,j} (d_{i,j} )\left(g_{2} (v_{j} \sin (\theta _{j} ))-g_{2} (v_{i} \sin (\theta _{i} ))\right)  
\end{aligned}\label{GrindEQ__3_13_} \\
&\begin{aligned}
\Lambda _{i} (w)=&\sum _{j\ne i}V'_{i,j}  (d_{i,j} )\frac{(x_{i} -x_{j} )}{d_{i,j} }\\
& -\sum _{j\ne i}\kappa _{i,j} (d_{i,j} )\left(g_{1} (v_{j} \cos (\theta _{j} ))-g_{1} (v_{i} \cos (\theta _{i} ))\right)  
 \end{aligned}\label{GrindEQ__3_14_}\end{align}\begin{align} 
&k_{i} (w)=\mu _{2} +\frac{\Lambda _{i} (w)}{v^{*} } +\frac{v_{\max } \cos (\theta _{i} )}{v^{*} (v_{\max } \cos (\theta _{i} )-v^{*} )} r\left(-\Lambda _{i} (w)\right) \label{GrindEQ__3_15_} 
\end{align}
and $\mu _{1} ,\mu _{2} >0$ are constants (controller gains), and $r\in C^{1} (\mathbb{R} )$, $g_{j} \in C^{1} (\mathbb{R} )$, $j=1,2$, are functions that satisfy
\begin{align}
\max (0,x)\le& r(x) \;\textrm{ for all }x\in \mathbb{R} \label{GrindEQ__3_16_}\\
g'_{j} (x)>&0 \;\;\; {\rm for\; }\, x\in \mathbb{R} ,\, \, j=1,2.   \label{GrindEQ__3_17_} 
\end{align}

The term $k_{i} (w)$ in the acceleration $F_{i} (t)$ given by \eqref{GrindEQ__3_15_}, is a state-dependent controller gain, which guarantees that the speed of each vehicle will remain positive and less than the speed limit. The first term in \eqref{GrindEQ__3_12_} drives the longitudinal speed of a vehicle towards the speed set-point $v^{*} $.  If $V_{i,j} $ in \eqref{GrindEQ__3_1_}, \eqref{GrindEQ__3_2_} is monotone, then, if vehicle $j$ is behind vehicle $i$, i.e., $(x_{i} -x_{j} )>0$, we have that $-V'_{i,j} (d_{i,j} )\frac{(x_{i} -x_{j} )}{d_{i,j} } >0$, and this term represents the effect of nudging \cite{19}, since vehicles that are close and behind vehicle $i$ will exert on it a ``pushing'' force that increases its acceleration. Notice that the fact that $b>1-\frac{v^{*}}{v_{\max}}$ implies that the denominator in \eqref{GrindEQ__3_11_} is positive for all $v_i\in(0,v_{\max})$, $\theta_i\in(-\varphi,\varphi)$, $i=1,...,n$. It should be noted that properties \eqref{GrindEQ__3_2_} and \eqref{GrindEQ__3_6_} guarantee that the feedback laws \eqref{GrindEQ__3_11_} and \eqref{GrindEQ__3_12_}, are \textit{decentralized} (per vehicle) and depend only on the relative speeds and distances from adjacent vehicles, more precisely from vehicles that are located at a distance less than $\lambda >0$.  

When the NCC is inviscid, then it does not require measurement of the speeds of the adjacent vehicles.

\subsection{Pseudo-Relativistic Cruise Controller (PRCC)}

The CLF in this case is given by 
\begin{equation} \label{GrindEQ__3_18_} 
\begin{aligned}{H_{R} (w):} & {=}  \frac{1}{2} \sum _{i=1}^{n}\frac{\left(v_{i} \cos (\theta _{i} )-v^{*} \right)^{2} +bv_{i}^{2} \sin ^{2} (\theta _{i} )}{\left(v_{\max } -v_{i} \right)v_{i} } \\ &+\sum_{i=1}^nU_i(y_i)+\frac{1}{2}\sum_{i=1}^n\sum_{i\neq j}V_{i,j}(d_{i,j})\\
&+A\sum_{i=1}^n\left(\frac{1}{cos(\theta_i)-\cos(\varphi)}-\frac{1}{1-\cos{\varphi}}\right) \end{aligned} 
\end{equation} 
where $A>0$, $b>1-\frac{v^{*} }{v_{\max } } >0$ (recall that $v^{*} \in (0,v_{\max } )$) are parameters of the controller and the Lyapunov function.  Notice that the kinetic energy term in $H_{R} $ (i.e., the term $\frac{1}{2} \sum _{i=1}^{n}\frac{\left(v_{i} \cos (\theta _{i} )-v^{*} \right)^{2} +bv_{i}^{2} \sin ^{2} (\theta _{i} )}{\left(v_{\max } -v_{i} \right)v_{i} }  $) is not related to the kinetic energy of classical mechanics, but is similar to the kinetic energy of a system of $n$ particles in relativistic mechanics, with speed limits 0 and $v_{\max } $ in place of $-c$ and $c$, where $c$ is the speed of light, which are the speed limits in relativistic mechanics. In relativistic mechanics, the kinetic energy increases to infinity when the speed of an object approaches (in absolute value) the speed of light, which indicates that no object with mass can reach the speed of light. Analogously, in \eqref{GrindEQ__3_18_}, the kinetic energy term grows to infinity as the speed of a vehicle approaches zero or the maximum speed $v_{\max } $, thus restricting the speed of vehicles in $(0,v_{\max } )$. As in the case of \eqref{GrindEQ__3_9_}, the sum of the terms $\sum _{i=1}^{n}U_{i} (y_{i} ) +\frac{1}{2} \sum _{i=1}^{n}\sum _{j\ne i}V_{i,j} (d_{i,j} )  $ is related to the potential energy of the system; and the last term of \eqref{GrindEQ__3_18_} ($A\sum _{i=1}^{n}\left(\frac{1}{\cos (\theta _{i} )-\cos (\varphi )} -\frac{1}{1-\cos (\varphi )} \right) $) is a penalty term that blows up when $\theta _{i} \to \pm \varphi $. 

The following result shows that the function $H_{R} $ is a size function for the state space $\Omega$ defined by \eqref{GrindEQ__2_7_}.

\noindent \textbf{Lemma 2: }\textit{Let constants $A>0$, $v_{\max } >0$, $v^{*} \in (0,v_{\max } )$,  $L_{i,j} >0$, $i,j=1,...,n$, $i\ne j$, $\lambda >0$ that satisfies \eqref{GrindEQ__3_8_}, $\varphi \in \left(0,\frac{\pi }{2} \right)$ that satisfies \eqref{GrindEQ__2_4_}, and define the function $H_{R} :\Omega \to \mathbb{R} _{+} $ by means of \eqref{GrindEQ__3_18_}, where $\Omega $ is given by \eqref{GrindEQ__2_7_}. Then, there exist non-decreasing functions $\eta _{i} :\mathbb{R} _{+} \to [0,a)$, $i=1,...,n$, }$\ell _{2} :\mathbb{R} _{+} \to [v^{*} ,v_{\max } )$, \textit{$\omega :\mathbb{R} _{+} \to \left[0,\varphi \right)$, a non-increasing function $\ell _{1} :\mathbb{R} _{+} \to (0,v^{*} ]$,} \textit{and, for each pair $i,j=1,...,n$, $i\ne j$, there exist non-increasing functions $\rho _{i,j} :\mathbb{R} _{+} \to (L_{i,j} ,\lambda ]$ with $\rho _{i,j} (s)\equiv \rho _{j,i} (s)$, such that the following implications hold for all $i,j=1,...,n$, $ j\ne i$:}
\begin{align}
w\in \Omega \Rightarrow &\ell _{1} \left(H_{R} (w)\right)\le v_{i} \le \ell _{2} \left(H_{R} (w)\right)\, ,\left|\theta _{i} \right|\le \omega \left(H_{R} (w)\right),\nonumber\\
&\left|y_{i} \right|\le \eta _{i} \left(H_{R} (w)\right),d_{i,j} \ge \rho _{i,j} \left(H_{R} (w)\right),\label{GrindEQ__3_19_}
\end{align}

\noindent Let $f_{j} :\mathbb{R} \to \mathbb{R} $, $j=1,2$, be $C^{1} $ functions that satisfy:
\begin{equation}
f_{j} (0)=0 \textrm{ and  }x\, f_{j} (x)>0,\, \, {\rm for\; }\, x\ne 0,\, \, j=1,2  \label{GrindEQ__3_20_}
\end{equation}
and $g_{j} :\mathbb{R} \to \mathbb{R} $, $j=1,2$, be $C^{1} $ functions that satisfy \eqref{GrindEQ__3_17_}. The controllers that correspond to the CLF \eqref{GrindEQ__3_18_} are
\begin{align}
&\begin{aligned}
u_{i} =&\frac{v_{i} }{\beta (v_{i} ,\theta _{i} )} \Biggl(G_{i} (w)-U'_{i} (y_{i} )-a(v_{i} ,\theta _{i} )F_{i} \\ & -\sum _{j\ne i}p_{i,j} V'_{i,j} (d_{i,j} )\frac{\left(y_{i} -y_{j} \right)}{d_{i,j} }  \Biggr) 
\end{aligned} \label{GrindEQ__3_21_} \\
&F_{i} =\frac{1}{q(v_{i} ,\theta _{i} )} \left(R_{i} (w)-\sum _{j\ne i}V'_{i,j} (d_{i,j} )\frac{\left(x_{i} -x_{j} \right)}{d_{i,j} }  \right) \label{GrindEQ__3_22_} 
\end{align} 
where 
\begin{align}
\begin{aligned}
G_{i}& (w)=-f_{2} \left(v_{i} \sin (\theta _{i} )\right)\\&+\sum _{j\ne i}\kappa _{i,j} \left(d_{i,j} \right)\left(g_{2} \left(v_{j} \sin (\theta _{j} )\right)-g_{2} \left(v_{i} \sin (\theta _{i} )\right)\right)  
\end{aligned} \label{GrindEQ__3_23_} \\
\begin{aligned}
R_{i}& (w)=-f_{1} \left(v_{i} \cos (\theta _{i} )-v^{*} \right)\\&+\sum _{j\ne i}\kappa _{i,j} \left(d_{i,j} \right)\left(g_{1} \left(v_{j} \cos (\theta _{j} )\right)-g_{1} \left(v_{i} \cos (\theta _{i} )\right)\right)  
\end{aligned}\label{GrindEQ__3_24_} 
\end{align}
and 
\begin{align}
q(v,\theta ):=&\frac{v_{\max } v\cos (\theta )+v^{*} v_{\max } -2v^{*} v}{2\left(v_{\max } -v\right)^{2} v^{2} }  \label{GrindEQ__3_25_}  \\
\beta (v,\theta ):=&\frac{A}{\left(\cos (\theta )-\cos (\varphi )\right)^{2} } +\frac{(b-1)v\cos (\theta )+v^{*} }{\left(v_{\max } -v\right)}  \label{GrindEQ__3_26_} \\
a(v,\theta ):=&\frac{bv_{\max } \sin (\theta )}{2\left(v_{\max } -v\right)^{2} v}  .  \label{GrindEQ__3_27_} 
\end{align} 
Notice that the definition of $b$ implies that $\beta (v,\theta )>0$ for all $v\in (0,v_{\max } )$, and $\theta \in (-\varphi ,\varphi )$. 

The pseudo-relativistic feedback laws \eqref{GrindEQ__3_21_} and \eqref{GrindEQ__3_22_} are derived by using the CLF \eqref{GrindEQ__3_18_}, which is also a size function.  Compared to the Newtonian controller \eqref{GrindEQ__3_11_}, \eqref{GrindEQ__3_12_}, the controller \eqref{GrindEQ__3_21_}, \eqref{GrindEQ__3_22_} is simpler, since it does not use state-dependent controller gains to restrict the speed in $(0,v_{\max } )$ (due to the properties of the size function $H_{R} $). The function $\frac{1}{q(v,\theta )} $ in \eqref{GrindEQ__3_22_}, \eqref{GrindEQ__3_25_} drives the acceleration $F_{i} $ to zero, when the speed of the vehicle tends to zero or to the maximum speed $v_{\max } $.   Finally, notice that properties \eqref{GrindEQ__3_2_} and \eqref{GrindEQ__3_6_} guarantee that the feedback laws \eqref{GrindEQ__3_21_} and \eqref{GrindEQ__3_22_}, are \textit{decentralized} (per vehicle) and depend only on the relative speed of and distance from adjacent vehicles, namely vehicles that are located at a distance less than $\lambda >0$.  

Again, it should be noted that, when the PRCC is inviscid, then it does not require measurement of the speeds of the adjacent vehicles.

\subsection{Main Results}

The following theorem shows that each of the closed-loop systems \eqref{GrindEQ__2_2_}, \eqref{GrindEQ__3_11_}, \eqref{GrindEQ__3_12_}, and \eqref{GrindEQ__2_2_}, \eqref{GrindEQ__3_21_}, \eqref{GrindEQ__3_22_} satisfy properties (P1) and (P2).

\noindent \textbf{Theorem 1 (Closed-loop system with PRCC):} \textit{For every $w_{0} \in \Omega $ there exists a unique solution $w(t)\in \Omega $ of the initial-value problem \eqref{GrindEQ__2_2_}, \eqref{GrindEQ__3_21_}, \eqref{GrindEQ__3_22_} with initial condition $w(0)=w_{0} $. The solution $w(t)\in \Omega $ is defined for all $t\ge 0$ and satisfies, for $i=1,...,n$,}
\begin{align}
\mathop{\lim }\limits_{t\to +\infty } \left(v_{i} (t)\right)=v^{*} , \mathop{\lim }\limits_{t\to +\infty } \left(\theta _{i} (t)\right)=0 \label{GrindEQ__3_28_} \\
\mathop{\lim }\limits_{t\to +\infty } \left(u_{i} (t)\right)=0, \mathop{\lim }\limits_{t\to +\infty } \left(F_{i} (t)\right)=0.             \label{GrindEQ__3_29_}                                     
\end{align} 
\textit{Moreover, there exist non-decreasing functions $Q_{k} :\mathbb{R} _{+} \to \mathbb{R} _{+} $ ($k=1,2$) such that the following inequalities hold for every solution $w(t)\in \Omega $ of \eqref{GrindEQ__2_2_}, \eqref{GrindEQ__3_21_}, \eqref{GrindEQ__3_22_}}
\begin{align}
\mathop{\max }\limits_{i=1,...,n} \left(\left|F_{i} (t)\right|\right)\le Q_{1} (H_{R} (w(0))), \textrm{ for all }t\ge 0                                       \label{GrindEQ__3_30_}\\
\mathop{\max }\limits_{i=1,...,n} \left(\left|u_{i} (t)\right|\right)\le Q_{2} (H_{R} (w(0))), \textrm{ for all }t\ge 0                                       \label{GrindEQ__3_31_}
\end{align}

\noindent \textbf{Theorem 2 (Closed-loop system with NCC):} \textit{For every $w_{0} \in \Omega $ there exists a unique solution $w(t)\in \Omega $ of the initial-value problem \eqref{GrindEQ__2_2_}, \eqref{GrindEQ__3_11_}, \eqref{GrindEQ__3_12_} with initial condition $w(0)=w_{0} $. The solution $w(t)\in \Omega $ is defined for all $t\ge 0$ and satisfies \eqref{GrindEQ__3_28_} and \eqref{GrindEQ__3_29_} for all $i=1,...,n$. Moreover, there exist non-decreasing functions $Q_{k} :\mathbb{R} _{+} \to \mathbb{R} _{+} $ ($k=3,4$) such that the following inequalities hold for every solution $w(t)\in \Omega $ of \eqref{GrindEQ__2_2_}, \eqref{GrindEQ__3_11_}, \eqref{GrindEQ__3_12_}:}
\begin{align}
\mathop{\max }\limits_{i=1,...,n} \left(\left|F_{i} (t)\right|\right)\le Q_{3} (H(w(0))), \textrm{ for all }t\ge 0                                       \label{GrindEQ__3_32_}\\
\mathop{\max }\limits_{i=1,...,n} \left(\left|u_{i} (t)\right|\right)\le Q_{4} (H(w(0))), \textrm{ for all }t\ge 0                                       \label{GrindEQ__3_33_}
\end{align}

\textbf{Remarks:} 
\begin{enumerate}
\item  The results of Theorem 1 and Theorem 2 hold globally, i.e., for any initial condition $w_{0} \in \Omega $.

\item  It is important to notice that, due to technical constraints, an inequality of the form $\left|F_{i} (t)\right|\le K$ must be satisfied for all $t\ge 0$, where $K>0$ is a constant that depends on the technical characteristics of the vehicles and the road, as well as passenger convenience. Inequalities \eqref{GrindEQ__3_30_}, \eqref{GrindEQ__3_32_} allow us to determine the set of initial conditions $w_{0} \in \Omega $, for which the inequality $\left|F_{i} (t)\right|\le K$ holds: it includes the set of all $w_{0} \in \Omega $ with $Q_{1} (H_{R} (w_{0} ))\le K$ in the case of PRCC and $Q_{3} (H(w_{0} ))\le K$ in the case of NCC.  
\end{enumerate}

\subsection{Differences Between the Two Families of Controllers}

To illustrate the behaviour of the closed-loop systems with the NCCs and PRCCs, we consider the case of $n=10$ vehicles of equal length and we set $L_{i,j} =L$ and $p_{i,j} =p$ for all $i,j=1,...,n$. Let
\begin{align} 
&V_{i,j}(d)=\left\{\begin{array}{cc} {q_{1} \frac{(\lambda -d)^{3} }{d-L} } & {,L<d\le \lambda } \\ {0} & {,d>\lambda } \end{array}\right. \label{GrindEQ__3_34_} \\ 
&U_i(y)=\left\{\begin{array}{cc} {\left(\frac{1}{a^{2} -y^{2} } -\frac{c}{a^{2} } \right)^{4} } & {\begin{array}{l} {,-a<y<-\frac{a\sqrt{c-1} }{\sqrt{c} } \; \, \, } \\ {and\, \, \, \, \frac{a\sqrt{c-1} }{\sqrt{c} } <y<a\; } \end{array}} \\ {0} & {-\frac{a\sqrt{c-1} }{\sqrt{c} } \; \le y\le \, \, \, \frac{a\sqrt{c-1} }{\sqrt{c} } } \end{array}\right. \, \, \,   \label{GrindEQ__3_35_} \\
&\kappa_{i,j}(d)=\left\{\begin{array}{cc} {q_{2} (\lambda -d)^{2} } & {L<d\le \lambda } \\ {0} & {d>\lambda } \end{array}\right.   \label{GrindEQ__3_36_} 
\end{align} 
where $\lambda >L>0$, $c{\kern 1pt} \ge 1$ and $q_{1} >0,q_{2} \ge 0$. 

For the PRCC, we considered the case where $f_{1} (s)=\mu _{2} s$, $f_2(s)=\mu_1s$, $g_{k} (s)=s$, for $k=1,2$ with $\mu _{2} >0$ and $g_{k} (s)=s$, $k=1,2$. For the NCC, the function $r(x)$ satisfying \eqref{GrindEQ__3_15_} was selected to be 
\begin{equation} \label{GrindEQ__3_37_} 
r(x)=\frac{1}{2\varepsilon } \left\{\begin{array}{c} {\quad \quad \quad 0\, \, \quad \quad \quad \; \; \; if\; x\le -\varepsilon } \\ {\quad \quad \quad \left(x+\varepsilon \right)^{2} \; \quad \quad \quad if\; -\varepsilon <x<0} \\ {\varepsilon ^{2} +2\varepsilon x\; \quad \; \quad \quad if\; x\ge 0} \end{array}\right.  
\end{equation} 
where $\varepsilon >0$, $c{\kern 1pt} \ge 1$. 

Notice that, when $v_{i} \approx v^{*} $, $\theta _{i} \approx 0$ and $d\approx \lambda $, then $q(v,\theta )\approx v_{\max }^{2} $ and $k_{i} (w)\approx \mu _{2} $. In order to compare the two controllers, the gains of the controllers are set approximately equal. Thus, we selected $\mu _{2} =v_{\max }^{-2} $, $A=1,\, \, \beta =1$ and $q_{1} =v_{\max }^{-2} 10^{-3} $ for the PRCC, with $q_{2} =0.5v_{\max }^{-2} $ for the viscous case and $q_{2} =0$ for the inviscid case; and $q_{1} =10^{-3} $,$\mu _{2} =v_{\max }^{-1} $, $\mu _{1} =0.4$, $A=1,\, \, \beta =1$ for the NCC, with $q_{2} =0.5$ for the viscous case and $q_{2} =0$ for the inviscid case. Finally, we selected $c=1.5$, $\lambda =25$, $\varepsilon =0.2$, $v_{\max } =35$, $a=7.2$, and $\varphi =0.25$, $p=5.11$  and $L=5.59$.   

As a measure for comparison, we use the CLF $H$ defined by \eqref{GrindEQ__3_9_} for both the NCC and PRCC. Figure 1 shows the asymptotic convergence of $H$ for the inviscid case (solid line) and the viscous case (dashed line) for both the NCC and PRCC. It is shown that convergence of $H$ is faster in the viscous case than in the inviscid case for both cruise controllers. \textit{{This indicates that the additional measurement requirements of the viscous cruise controllers provide a direct advantage over the inviscid case, where no measurement of the speeds of the adjacent vehicles is required.}} We also checked the logarithm of $H$ and we found that convergence is not exponential for all controllers (viscous or inviscid). 

\noindent \includegraphics*[width=3.40in, height=1.89in ]{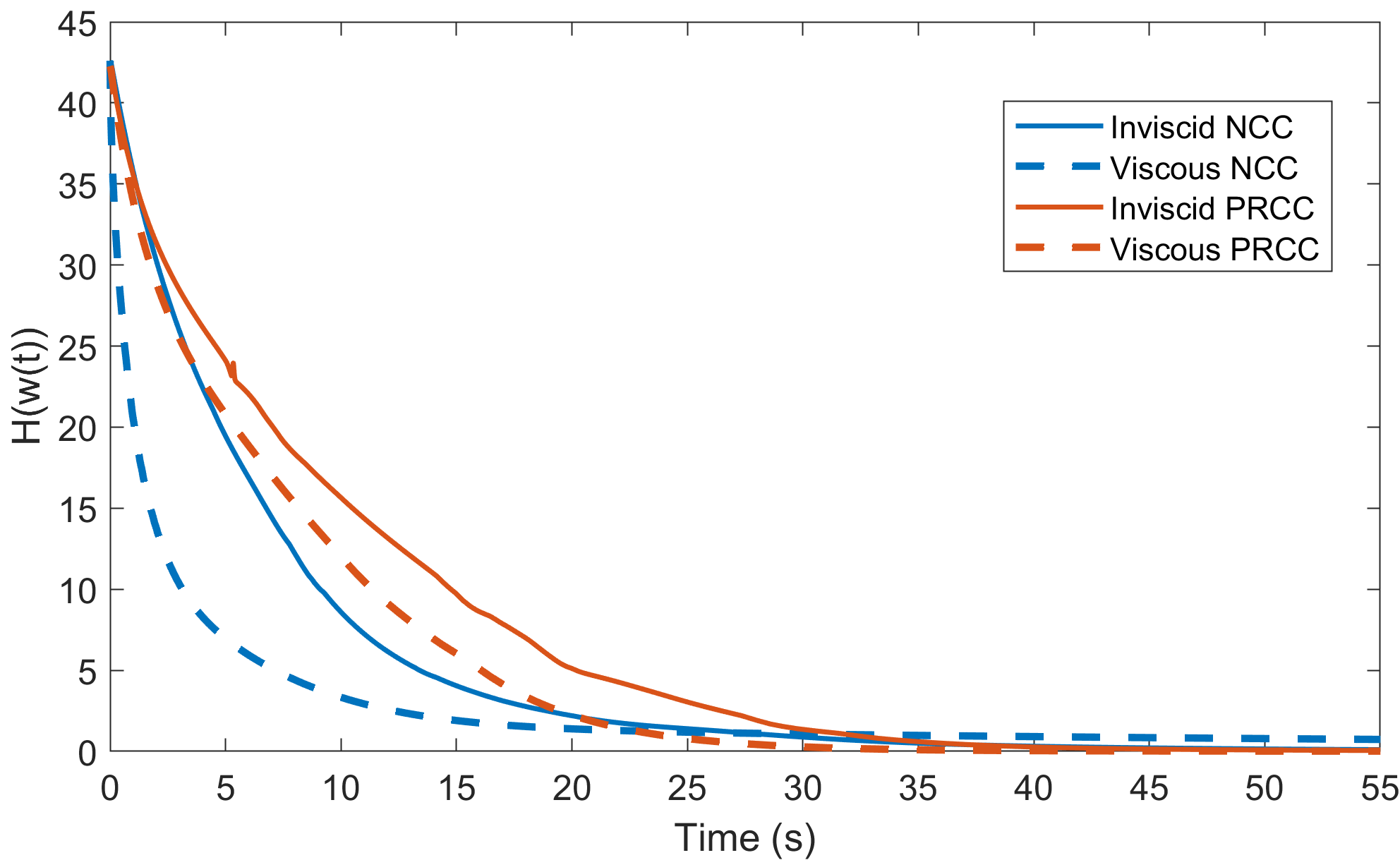}

\noindent \textbf{Figure 1: }Convergence of $H$ for the inviscid and viscous NCC and PRCC.

\section{Macroscopic Models}

In this section we present the macroscopic models that correspond to the microscopic model \eqref{GrindEQ__2_2_} with the NCC \eqref{GrindEQ__3_11_}, \eqref{GrindEQ__3_12_} or the PRCC \eqref{GrindEQ__3_21_}, \eqref{GrindEQ__3_22_}.  

Let $\rho _{\max } ,v_{\max } >0$ and $v^{*} \in (0,v_{\max } )$, $\bar{\rho }\in (0,\rho _{\max } )$ be constants and let $\mu :(0,\rho _{\max } )\to \mathbb{R} _{+} $,  $P:(0,\rho _{\max } )\to \mathbb{R} _{+} $ be $C^{2} \left((0,\rho _{\max } )\right)$ and non-negative functions that satisfy: 
\begin{align}
&\mathop{\lim \sup }\limits_{\rho \to \rho _{\max }^{-} } \left(P(\rho )\right)=+\infty ,   \label{GrindEQ__4_1_} \\
\mu (\rho )=&0, \,\,P(\rho )=0 \textrm{ for all }\rho \in \left(0,\bar{\rho }\right]  \label{GrindEQ__4_2_}
\end{align}

\noindent The macroscopic model that corresponds to the microscopic model \eqref{GrindEQ__2_2_} under the PRCC \eqref{GrindEQ__3_21_}, \eqref{GrindEQ__3_22_} is the following system of PDEs for $t>0$ and $x\in I(t)$, where $I(t)\subseteq \mathbb{R} $ is an appropriate interval 
\begin{align} 
\frac{\partial \, \rho }{\partial \, t} &+\frac{\partial }{\partial \, x} \left(\rho v\right)=0,   \label{GrindEQ__4_3_} \\
q(v)\frac{\partial \, v}{\partial \, t} +&q(v)v\frac{\partial \, v}{\partial \, x} +\frac{P'(\rho )}{\rho } \frac{\partial \, \rho }{\partial \, x} = \frac{1}{\rho } \frac{\partial \, }{\partial \, x} \left(\frac{\partial \, v}{\partial \, x} g'(v)\mu \left(\rho \right)\right)\nonumber\\&\qquad\qquad\qquad\qquad -f\left(v-v^{*} \right) 
\label{GrindEQ__4_4_}
\end{align} 
with constraints $0<\rho (t,x)<\rho _{\max } $, $0<v(t,x)<v_{\max } $ for all $t>0$ and $x\in I(t)$, where $g\in C^{1} \left(\mathbb{R}\right)$ is an increasing function with $g'(v)>0$ for all $v\in\mathbb{R}$, $f:\mathbb{R} \to \mathbb{R} $ with $f(0)=0$ is a $C^{1} $ function with $x\, f(x)>0$ for all $x\ne 0$ and 
\begin{equation} \label{GrindEQ__4_5_} 
q(v):=\frac{v_{\max } v-2v^{*} v+v^{*} v_{\max } }{2\left(v_{\max } -v\right)^{2} v^{2} } .   
\end{equation} 
The term $f(v-v^{*} )$ is a \textit{relaxation term }that describes the tendency of vehicles to adjust their speed to the speed set-point $v^{*} $ and is similar to friction, see \cite{6}. The term $\frac{P'(\rho )}{\rho } \frac{\partial \, \rho }{\partial \, x} $ is a \textit{pressure term} and expresses the tendency to accelerate or to decelerate based on the (local) density. Finally, the term $\frac{1}{\rho } \frac{\partial \, }{\partial \, x} \left(\frac{\partial \, v}{\partial \, x} g'(v)\mu \left(\rho \right)\right)$, is a \textit{viscosity term,} by analogy with the theory of fluids, with $\mu (\rho )$ playing the role of dynamic viscosity. When $g(v)\equiv v$, the viscosity term is exactly the same as the viscosity term appearing in compressible fluid flow (see \cite{13}, \cite{16}, \cite{24} and references therein). When $g(v)$ does not coincide with $v$, then the viscosity term is similar to the viscosity term appearing in porous fluid flow (see \cite{1}).

The macroscopic model that corresponds to the NCC is given by the continuity equation \eqref{GrindEQ__4_3_} and the following momentum equation for $t>0$ and $x\in I(t)$, where $I(t)\subseteq \mathbb{R} $ is an appropriate interval 
\begin{equation} \label{GrindEQ__4_6_} 
\begin{aligned}
\frac{\partial \, v}{\partial \, t} +v\frac{\partial \, v}{\partial \, x} +\frac{P'(\rho )}{\rho } \frac{\partial \rho }{\partial x} =&\frac{1}{\rho } \frac{\partial }{\partial x} \left(\frac{\partial v}{\partial x} g'(v)\mu (\rho )\right)\\
&-\left(\gamma +h\left(G\right)\right)\left(v-v^{*} \right) 
\end{aligned}
\end{equation} 
where
\begin{align*}
G=&-\frac{P'(\rho )}{\rho } \frac{\partial \rho }{\partial x} +\frac{1}{\rho } \frac{\partial }{\partial x} \left(\frac{\partial v}{\partial x} g'(v)\mu (\rho )\right)\\
h(s)=&\frac{v_{\max } r(s)}{v^{*} (v_{\max } -v^{*} )} -\frac{s}{v^{*} } 
\end{align*}
and $r\in C^{1} $ satisfies $\max (0,x)\le r(x)$ for all $x\in \mathbb{R} $. Again, the macroscopic model \eqref{GrindEQ__4_3_}, \eqref{GrindEQ__4_6_} is to be considered with constraints $0<\rho (t,x)<\rho _{\max } $, $0<v(t,x)<v_{\max } $  for all  $t>0$ and $x\in I(t)$, with $g\in C^{1} \left(\mathbb{R}\right)$ being an increasing function with $g'(v)>0$ for all $v\in \mathbb{R}$ and $\gamma >0$ being a constant. 

\textbf{    Remarks: (i) }Both models \eqref{GrindEQ__4_3_}, \eqref{GrindEQ__4_4_} and \eqref{GrindEQ__4_3_}, \eqref{GrindEQ__4_6_} do not include non-local terms and have certain characteristics from the kinematic theory of fluids. Traffic flow is isotropic, as in fluid flow, since the vehicles are autonomous and do not react based on downstream vehicles only (as in conventional traffic). Due to the nudging induced by the NCCs and PRCCs, vehicles are affected by both upstream and downstream vehicles. 

\noindent \textbf{(ii)} There are infinite equilibrium points for both models, namely the points where $v(x)\equiv v^{*} $ and $\rho (x)\le \bar{\rho }$ for all $x\in \mathbb{R} $.  

\noindent \textbf{(iii) }For the inviscid NCC-based model ($\mu (s)\equiv 0$), it was shown in \cite{12}, that, if the density is sufficiently small, then the solution of the macroscopic model approaches the equilibrium speed (in the sup norm); while the density converges exponentially to a traveling wave.

\noindent \textbf{(iv) }the model \eqref{GrindEQ__4_6_} is highly nonlinear due to the presence of a highly nonlinear relaxation term $\left(\gamma +h\left(G\right)\right)\left(v-v^{*} \right)$ in the speed PDE.


The macroscopic models \eqref{GrindEQ__4_3_}, \eqref{GrindEQ__4_4_} and \eqref{GrindEQ__4_3_}, \eqref{GrindEQ__4_6_} can approximate the movement of $n$ identical vehicles with total mass $m>0$ on a straight road under the PRCC \eqref{GrindEQ__3_21_}, \eqref{GrindEQ__3_22_} and under the NCC \eqref{GrindEQ__3_11_}, \eqref{GrindEQ__3_12_} when the following assumption holds: 

\textbf{Assumption H}: (i) the vehicles are constrained to move on a line (longitudinal motion), 

\noindent (ii) there exist constants $\lambda >L>0$ with $\lambda <2L$ such that $V_{i,j} (s)=\Phi \left(ns\right)$ for all $i,j=1,...,n$ and $s>L/n$, where $\Phi :(L,+\infty )\to \mathbb{R} _{+} $ is a $C^{2} $ function that satisfies $\mathop{\lim }\limits_{d\to L^{+} } \left(\Phi (d)\right)=+\infty $ and $\Phi (d)=0$ for all $d\ge \lambda $, 

\noindent (iii) $\kappa _{i,j} (s)=n^{2} K(ns)$, for all $i,j=1,...,n$ and $s>L/n$, where $K:(L,+\infty )\to \mathbb{R} _{+} $ is a $C^{1} $ function that satisfies $K(d)=0$ for all $d\ge \lambda $, and 

\noindent (iv) the number of vehicles $n$ is very large (tends to infinity).  

Note that assumption (i) above, is not always accurate since it neglects the lateral movement of vehicles.

We next present the relations between the various parameters and functions involved in the NCC and PRCC on one hand; and the corresponding macroscopic quantities involved in the macroscopic models on the other hand. Table 1 shows how all parameters and functions of the macroscopic models can be directly obtained from the corresponding microscopic models.

Notice that the pressure $P$ is based on the derivative of the potential $\Phi $, which, in the microscopic case, exerts the same force to the following and preceding vehicle (nudging). The latter is in analogy with fluids where pressure at any point in a fluid is (locally) the same in all directions. The dynamic viscosity $\mu (\rho )$ makes the ``traffic fluid'' act as a Newtonian fluid. However, in contrast to actual fluids, the ``traffic fluid'' induced by the PRCCs or the NCCs satisfies \eqref{GrindEQ__4_2_}. The function $g(v)$ is met in various nonlinear PDEs. For fluid flows in non-porous media, we have $g(v)\equiv v$. For compressible fluid flows in porous media, the function $g(v)$ takes the form $g(v)=cv^{m} $, where $c,m>0$ are constants (see \cite{1}, but notice also the essential difference that, for flow in porous media, the function $g$ is a function of the density $\rho $ instead of the speed $v$). The relaxation terms appear as friction terms (see \cite{6}, but, instead of penalizing the speed $v$, the relaxation terms penalize the deviation of the speed from the desired speed $v^{*} $).  

\begin{table}
\centering
\begin{tabular}{|p{1.15in}|c|p{1.0in}|} \hline 
 & Macroscopic  & Microscopic  \\ \hline 
Maximum   Density & $\rho _{\max } $  & $=\frac{m}{L} $  \\ \hline 
Maximum Speed & $v_{max}$ & $=v_{\max } $ \\ \hline 
Desired Speed & $v^{*} $ & $=v^{*} $ \\ \hline 
  Interaction \newline   Density & $\bar{\rho }$  & $=\frac{m}{\lambda } $ \\ \hline 
Dynamic Viscosity  & $\mu (\rho )$  & $=\frac{m^{2} }{\rho } K\left(\frac{m}{\rho } \right)$ \\ \hline 
 & $g(v)$  & $=g_{1} (v)$  \\ \hline 
Pressure & $P(\rho )$  & $=z-m\Phi '\left(\frac{m}{\rho } \right)$, \newline arbitrary $z\in \mathbb{R} $   \\ \hline 
Relaxation Term for the PRCC & $f(v-v^{*} )$  & $=f_{1} (v-v^{*} )$  \\ \hline 
Relaxation Term for the NCC & $r(x)$  & $=r(x)$  \\ \hline 
Constant Parameter for the NCC & $\gamma$  & $=\mu _{2} $  \\ \hline 
\end{tabular}
\caption{Relations between macroscopic and microscopic quantities}
\end{table}
The most important implication from Table 1 is the fact that, by changing the functions and the parameters of the NCCs and the PRCCs, \textit{{we can directly determine the physical properties of the ``traffic fluid''.}} In this sense, we may talk about an engineered or designed artificial fluid that approximates the actual emerging traffic flow. To understand how far the implications of the relations between the cruise controller parameters and the characteristics of the traffic fluid go, it is important to notice that for isentropic (or barotropic) flow of gases, the dynamic viscosity and the pressure are always increasing functions of the fluid density (see the discussion in \cite{13}). However, for a traffic fluid that emerges from the use of NCCs or PRCCs, if the cruise controller uses a non-monotone potential function $\Phi $ or a non-monotone viscosity function $K$, then it is possible to obtain a traffic fluid with dynamic viscosity and pressure, which are non-monotone functions of the fluid density (and can have local minima). Thus, the traffic fluid can be arranged to have very different physical properties from those of real compressible fluids (gases).    

Finally, an interesting case arises for a particular selection of $g$ satisfying $g'(v)>0$ for all $v\in \mathbb{R}$. We define the new speed-like variable 
\begin{equation} \label{GrindEQ__4_7_} 
w=g(v) 
\end{equation} 
where 
\begin{equation} \label{GrindEQ__4_8_} 
g(v)=\int _{v^{*} }^{v}q(s)ds .                                                   
\end{equation} 
The inverse function $g^{-1} $ is well-defined, since $g$ is increasing, hence, using \eqref{GrindEQ__4_7_}, \eqref{GrindEQ__4_8_}, the macroscopic model \eqref{GrindEQ__4_3_}, \eqref{GrindEQ__4_4_} can be written in the following form

\noindent 
\begin{align}
\frac{\partial \rho }{\partial t} +&\frac{\partial }{\partial x} \left(\rho g^{-1} (w)\right)=0  \label{GrindEQ__4_9_} \\
\rho \frac{\partial \, w}{\partial \, t} +\rho g^{-1} (w)\frac{\partial \, w}{\partial \, x} =&\frac{\partial }{\partial \, x} \left(\mu (\rho )\frac{\partial \, w}{\partial \, x} -P(\rho )\right)\nonumber\\&-\rho f\left(g^{-1} (w)\right).   \label{GrindEQ__4_10_} 
\end{align}

\noindent The above model \eqref{GrindEQ__4_9_}, \eqref{GrindEQ__4_10_} is very similar to isentropic compressible fluid flow models for Newtonian fluids (see \cite{13}, \cite{16}, \cite{24}) with a difference in the inertial term in \eqref{GrindEQ__4_10_} (where we have $\rho g^{-1} (w)\frac{\partial \, w}{\partial \, x} $ instead of $\rho w\frac{\partial \, w}{\partial \, x} $) and the flow term in \eqref{GrindEQ__4_9_} (where we have $\frac{\partial }{\partial x} \left(\rho g^{-1} (w)\right)$ instead of $\frac{\partial }{\partial x} \left(\rho w\right)$). The relaxation term $-\rho f\left(g^{-1} (w)\right)$ is a term similar to friction (see \cite{6}).

\section{Proofs}

The proof of Lemma 1 is a straightforward modification of Proposition 1 in \cite{11} and is omitted.

\textbf{Proof of Lemma 2: }The existence of the functions $\eta _{i} :\mathbb{R} _{+} \to [0,a)$, $\rho _{i,j} :\mathbb{R} _{+} \to (L_{i,j} ,\lambda ]$, $i,j=1,...,n$, $i\ne j$, and $\omega :\mathbb{R} _{+} \to \left[0,\varphi \right)$ for which implication \eqref{GrindEQ__3_19_} holds, is a direct consequence of Lemma 1 with $H_{R} (w)$ in place of $H(w)$. ${\rm I}$t suffices to show that the first inequality in \eqref{GrindEQ__3_19_} holds. Define the function $f(v):=\frac{(v-v^{*} )^{2} }{(v_{\max } -v)v} $ which is increasing on $[v^{*} ,v_{\max } )$ and decreasing on $(0,v^{*} ]$ with $f(v^{*} )=0$ and satisfies
\begin{equation}
f(v_{i} )\le H_{R} (w) \textrm{ for all } w\in \Omega.  \label{GrindEQ__5_6_}
\end{equation}

\noindent Define the function $\ell _{1} :\mathbb{R} _{+} \to (0,v^{*} ]$ to be equal to the inverse function of $f$ restricted on $(0,v^{*} ]$. The function $\ell _{1} :\mathbb{R} _{+} \to (0,v^{*} ]$ is decreasing with $\ell _{1} (0)=v^{*} $ and $\mathop{\lim }\limits_{v\to +\infty } \ell _{1} (v)=0$. Define also the function $\ell _{2} :\mathbb{R} _{+} \to [v^{*} ,v_{\max } )$ to be equal to the inverse function of $f$ restricted on $[v^{*} ,v_{\max } )$. The function $\ell _{2} :\mathbb{R} _{+} \to [v^{*} ,v_{\max } )$ is increasing with $\ell _{2} (0)=v^{*} $ and $\mathop{\lim }\limits_{v\to +\infty } \ell _{2} (v)=v_{\max } $. Then, due to definitions of $\ell _{1} ,\ell _{2} $ above and \eqref{GrindEQ__5_6_} it holds that 
\begin{equation} \label{GrindEQ__5_7_} 
0<\ell _{1} (H_{R} (w))\le v_{i} \le \ell _{2} (H_{R} (w))<v_{\max } .   
\end{equation} 
This completes the proof. $\triangleleft $ 

The proof of Theorem 1 is performed by using Barbalat's lemma (\cite{14}) and its following variant which uses uniform continuity of the derivative of a function. 

\noindent \textbf{     Lemma 3: }\textit{If a function $g\in C^{2} \left(\mathbb{R} _{+} \right)$ satisfies $\mathop{\lim }\limits_{t\to +\infty } (g(t))$ $\in \mathbb{R} $ and $\mathop{\sup }\limits_{t\ge 0} \left(\left|\ddot{g}(t)\right|\right)<+\infty $, then, $\mathop{\lim }\limits_{t\to +\infty } \left(\dot{g}(t)\right)=0$.}

\noindent \textbf{     Proof of Theorem 1:} Let $w_{0} \in \Omega $ and consider the unique solution $w(t)$ of the initial value problem \eqref{GrindEQ__2_2_}, \eqref{GrindEQ__3_21_}, \eqref{GrindEQ__3_22_}. Using the fact that the set $\Omega $ is open (recall definitions (2.3), \eqref{GrindEQ__2_7_}), we conclude that there exists $t_{\max } \in (0,+\infty ]$ such that the solution $w(t)$ of \eqref{GrindEQ__2_2_}, \eqref{GrindEQ__3_21_}, \eqref{GrindEQ__3_22_} is defined on $\left[0,t_{\max } \right)$ and satisfies $w(t)\in \Omega $ for all $t\in \left[0,t_{\max } \right)$. Furthermore, if $t_{\max } <+\infty $ then there exists an increasing sequence of times $\left\{\, t_{i} \in \left[0,t_{\max } \right)\, :\, i=1,2,...\right\}$ with $\mathop{\lim }\limits_{i\to +\infty } \left(t_{i} \right)=t_{\max } $ and either $\mathop{\lim }\limits_{i\to +\infty } \left(dist\left(w(t_{i} ),\partial \Omega \right)\right)=0$ or $\mathop{\lim }\limits_{i\to +\infty } \left(\left|w(t_{i} )\right|\right)=+\infty $.

We show first, that the solution $w(t)$ of the initial value problem \eqref{GrindEQ__2_2_}, \eqref{GrindEQ__3_21_}, \eqref{GrindEQ__3_22_} satisfies $w(t)\in \Omega $ for all $t\ge 0$. By using \eqref{GrindEQ__2_2_}, \eqref{GrindEQ__2_5_}, and \eqref{GrindEQ__3_18_}, it follows that for all $w\in \Omega $ 
\begin{equation} \label{GrindEQ__5_8_} 
\begin{aligned}
&\nabla H_{R} (w)\dot{w}=\\
&\sum _{i=1}^{n}F_{i} \frac{\left(v_{i} \cos (\theta _{i} )-v^{*} \right)}{2\left(v_{\max } -v_{i} \right)^{2} v_{i}^{2} } \left(v_{\max } v_{i} \cos (\theta _{i} )+v^{*} v_{\max } -2v^{*} v_{i} \right)  \\ 
&+\sum _{i=1}^{n}\left(U'_{i} (y_{i} )+\frac{bv_{\max } \sin (\theta _{i} )F_{i} }{2\left(v_{\max } -v_{i} \right)^{2} v_{i} } \right.\\
& \left.+\left(\frac{A}{\left(\cos (\theta _{i} )-\cos (\varphi )\right)^{2} } +\frac{(b-1)v_{i} \cos (\theta _{i} )+v^{*} }{\left(v_{\max } -v_{i} \right)} \right)\frac{u_{i} }{v_{i} }  \right)\\&\times v_{i} \sin (\theta _{i} )  \\ 
&{+\frac{1}{2} \sum _{i=1}^{n}\sum _{j\ne i}V'_{i,j} (d_{i,j} )\frac{\left(x_{i} -x_{j} \right)\left(v_{i} \cos (\theta _{i} )-v_{j} \cos (\theta _{j} )\right)}{d_{i,j} }   } \\
&{+\frac{1}{2} \sum _{i=1}^{n}\sum _{j\ne i}p_{i,j} V'_{i,j} (d_{i,j} )\frac{\left(y_{i} -y_{j} \right)\left(v_{i} \sin (\theta _{i} )-v_{j} \sin (\theta _{j} )\right)}{d_{i,j} }   } \end{aligned} 
\end{equation} 
Using \eqref{GrindEQ__3_25_}, \eqref{GrindEQ__3_26_}, \eqref{GrindEQ__3_27_}, the fact that $\sum _{i=1}^{n}\sum_{j\ne i}V'_{i,j} (d_{i,j} )$ $\frac{(x_{i} -x_{j} )}{d_{i,j} }=0$, and \eqref{GrindEQ__5_8_}, we get:
\begin{equation} \label{GrindEQ__5_9_} 
\begin{aligned}&\nabla H_{R} (w)\dot{w}=\sum _{i=1}^{n}\left(U'_{i} (y_{i} )+\beta (v_{i} ,\theta _{i} )\frac{u_{i} }{v_{i} }\right. \\&\left.+a(v_{i} ,\theta _{i} )F_{i} +\sum _{j\ne i}p_{i,j} V'_{i,j} (d_{i,j} )\frac{\left(y_{i} -y_{j} \right)}{d_{i,j} }  \right)v_{i} \sin (\theta _{i} ) \\ 
&+\sum _{i=1}^{n}\left(v_{i} \cos (\theta _{i} )-v^{*} \right)\\
&\times\left(F_{i} q(v_{i} ,\theta _{i} )+\sum _{j\ne i}V'_{i,j} (d_{i,j} )\frac{\left(x_{i} -x_{j} \right)}{d_{i,j} }  \right)
\end{aligned} 
\end{equation} 
From \eqref{GrindEQ__5_9_}, \eqref{GrindEQ__3_21_}, \eqref{GrindEQ__3_22_}, properties \eqref{GrindEQ__3_3_}, \eqref{GrindEQ__3_7_}, and definitions \eqref{GrindEQ__3_25_}, \eqref{GrindEQ__3_26_}, and \eqref{GrindEQ__3_27_}, we further get that 
\begin{equation} \label{GrindEQ__5_10_} 
\begin{aligned}
&{\nabla H_{R} (w)\dot{w}} {=} -\sum _{i=1}^{n}\left(v_{i} \cos (\theta _{i} )-v^{*} \right)f_{1} \left(v_{i} \cos (\theta _{i} )-v^{*} \right)\\
& -\sum _{i=1}^{n}v_{i} \sin (\theta _{i} )f_{2} \left(v_{i} \sin (\theta _{i} )\right)\\
&-\frac{1}{2} \sum _{i=1}^{n}\sum _{j\ne i}\kappa _{i,j} \left(d_{i,j} (t)\right)\left(v_{j} (t)\sin (\theta _{j} (t))-v_{i} (t)\sin (\theta _{i} (t))\right)\\&\times\left(g_{2} \left(v_{j} (t)\sin (\theta _{j} (t))\right)-g_{2} \left(v_{i} (t)\sin (\theta _{i} (t))\right)\right)   \\
& -\frac{1}{2} \sum _{i=1}^{n}\sum _{j\ne i}\kappa _{i,j} \left(d_{i,j} (t)\right)\left(v_{j} (t)\cos (\theta _{j} (t))-v_{i} (t)\cos (\theta _{i} (t))\right)\\&\left(g_{1} \left(v_{j} (t)\cos (\theta _{j} (t))\right)-g_{1} \left(v_{i} (t)\cos (\theta _{i} (t))\right)\right)      \end{aligned} 
\end{equation} 
It follows from \eqref{GrindEQ__3_17_}, \eqref{GrindEQ__3_20_} and \eqref{GrindEQ__5_10_} that for all $w\in \Omega $
\begin{equation} \label{GrindEQ__5_11_} 
\nabla H_{R} (w)\dot{w}\le 0 
\end{equation} 
Since $w(t)\in \Omega $ for all $t\in \left[0,t_{\max } \right)$, it follows from \eqref{GrindEQ__2_3_}, \eqref{GrindEQ__2_7_} that $v_{i} (t)\in \left(0,v_{\max } \right)$, $\theta _{i} (t)\in \left(-\varphi ,\varphi \right)$ for all $t\in \left[0,t_{\max } \right)$ and $i=1,...,n$. Thus, \eqref{GrindEQ__2_2_} implies that $0\le \dot{x}_{i} (t)\le v_{\max } $ for all $t\in \left[0,t_{\max } \right)$ and $i=1,...,n$. Moreover, inequality \eqref{GrindEQ__5_11_} implies that 
\begin{equation}
H_{R} (w(t))\le H_{R} (w_{0} ), \textrm{ for all }t\in \left[0,t_{\max } \right)\label{GrindEQ__5_12_}
\end{equation}
\noindent Consequently, we obtain from \eqref{GrindEQ__3_19_} that for all $t\in \left[0,t_{\max } \right)$ and $i,j=1,...,n$, $j\ne i$:
\begin{equation} \label{GrindEQ__5_13_} 
\begin{aligned}
&{\ell _{1} \left(H_{R} (w_{0} )\right)\le v_{i} (t)\le \ell _{2} \left(H_{R} (w_{0} )\right)\, ,\left|\theta _{i} (t)\right|\le \omega \left(H_{R} (w_{0} )\right),} \\ 
&\left|y_{i} (t)\right|\le \eta_i \left(H_{R} (w_{0} )\right),x_{i} (0)\le x_{i} (t)\le x_{i} (0)+v_{\max } t\; ,\;\\& d_{i,j} (t)\ge \rho _{i,j} \left(H_{R} (w_{0} )\right)>L_{i,j}  
\end{aligned} 
\end{equation} 
which imply that for every increasing sequence of times $\left\{t_{i} \in [0,t_{\max } ):i=1,2,...\right\}$ with $\mathop{\lim }\limits_{i\to +\infty } (t_{i} )=t_{\max } $ we cannot have $\mathop{\lim }\limits_{i\to +\infty } (dist(w(t_{i} ),\partial \Omega ))=0$ or $\mathop{\lim }\limits_{i\to +\infty } (|w_{i} (t)|)=+\infty $. Therefore, $t_{\max } =\infty $.

\noindent Next, define
\begin{equation} \label{GrindEQ__5_14_}
\begin{aligned}
&\Delta (t):=\sum _{i=1}^{n}\left(v_{i} (t)\cos (\theta _{i} (t))-v^{*} \right)f_{1} \left(v_{i} (t)\cos (\theta _{i} (t))-v^{*} \right) \\
&+\sum _{i=1}^{n}v_{i} (t)\sin (\theta _{i} (t))f_{2} \left(v_{i} (t)\sin (\theta _{i} (t))\right)\\ 
&+\frac{1}{2} \sum _{i=1}^{n}\sum _{j\ne i}\kappa _{i,j} \left(d_{i,j} (t)\right)\left(v_{j} (t)\sin (\theta _{j} (t))-v_{i} (t)\sin (\theta _{i} (t))\right)\\
&\times\left(g_{2} \left(v_{j} (t)\sin (\theta _{j} (t))\right)-g_{2} \left(v_{i} (t)\sin (\theta _{i} (t))\right)\right) \\ 
&+\frac{1}{2} \sum _{i=1}^{n}\sum _{j\ne i}\kappa _{i,j} \left(d_{i,j} (t)\right)\left(v_{j} (t)\cos (\theta _{j} (t))-v_{i} (t)\cos (\theta _{i} (t))\right)\\
&\times\left(g_{1} \left(v_{j} (t)\cos (\theta _{j} (t))\right)-g_{1} \left(v_{i} (t)\cos (\theta _{i} (t))\right)\right) \end{aligned}
\end{equation}

\noindent Notice that definition \eqref{GrindEQ__5_14_} and \eqref{GrindEQ__5_10_}, \eqref{GrindEQ__5_11_} implies that $\Delta (t)=-\frac{d}{d\, t} H_{R} (w(t))\ge 0$ for all $t\ge 0$. Therefore, since $H_{R} (w)\ge 0$ for all $w\in \Omega $, we obtain:
\begin{equation} \label{GrindEQ__5_15_} 
\int _{0}^{\infty }\Delta (t)dt \le H_{R} (w_{0} ).   
\end{equation} 
In order to show that the solution of system \eqref{GrindEQ__2_2_}, \eqref{GrindEQ__3_21_}, \eqref{GrindEQ__3_22_} satisfies \eqref{GrindEQ__3_28_}, it suffices to show that there exists a constant $\bar{M}>0$ such that
\begin{equation} \label{GrindEQ__5_16_} 
\left|\frac{d}{dt} \left(\Delta (t)\right)\right|\le \bar{M}.   
\end{equation} 
Indeed, if \eqref{GrindEQ__5_16_} holds for some constant $\bar{M}>0$,  then, $\Delta (t)$ is uniformly continuous and due to \eqref{GrindEQ__5_15_}, we can apply Barbalat's lemma (\cite{14}), to obtain that $\mathop{\lim }\limits_{t\to +\infty } \Delta (t)=0$. Then, by using assumption \eqref{GrindEQ__3_20_} and the fact that $\Delta (t)\ge 0$ for all $t\ge 0$, it can be shown that \eqref{GrindEQ__3_28_} holds. 

In order to prove \eqref{GrindEQ__5_16_}, we show first that \eqref{GrindEQ__3_30_} and \eqref{GrindEQ__3_31_} hold. Define
\begin{align}
&\begin{aligned}
 B_{i,j} (s):=\max \left\{|V'_{i,j} (d)|:s\le d\le \lambda \right\} \textrm{ for }s\in (L_{i,j} ,\lambda ]&\\,\,\,i,j=1,...,n,\, j\ne i&
\end{aligned} \label{GrindEQ__5_17_}\\
&\begin{aligned}
c_{i,j} (s):=\max \left\{\kappa _{i,j} (d):s\le d\le \lambda \right\}\textrm{ for }s\in (L_{i,j} ,\lambda ],&\\ \,\,i,j=1,...,n,\, j\ne i&
\end{aligned}\label{GrindEQ__5_18_}\\
&\zeta _{i} (s):=\max \left\{|U'_{i} (y):|y|\le s\right\}\textrm{ for }s\in [0,a).      \label{GrindEQ__5_19_}
\end{align}

Notice that definition \eqref{GrindEQ__2_5_} implies that 
\begin{equation}\label{GrindEQ__5_20_}
\left|\frac{x_{i} -x_{j} }{d_{i,j} } \right|\le 1\,\, \textrm{and }  \left|\frac{y_i-y_j}{d_{i,j}}\right|\leq\frac{1}{\sqrt{p_{i,j}}}\,\,
\end{equation}
for all $w\in\Omega$, $i,j=1,\ldots,n,j\neq i$. Moreover, for each $i=1,...,n$, let $m_{i} \ge 2$ be the maximum number of points that can be placed within the area bounded by two concentric ellipses with semi-major axes $\bar{L}_{i} =\min \{L_{i,j} ,j=1,...,n,j\ne i\}$ and $\lambda $ satisfying \eqref{GrindEQ__3_8_}, and semi-minor axes $\frac{\bar{L}_{i} }{\mathop{\max }\limits_{j\ne i} \sqrt{p_{i,j} } } $ and $\frac{\lambda }{\mathop{\min }\limits_{j\ne i} \sqrt{p_{i,j} } } $  so that each point has distance (in the metric given by (2.5)) at least $\bar{L}_{i} $ from every other point. Then, it follows from \eqref{GrindEQ__2_5_}, \eqref{GrindEQ__3_2_}, the fact that $d_{i,j} >L_{i,j} $ for $i,j=1,...,n$, $j\ne i$ and the definition of $m_{i} $ above, that the sums $\sum _{j\ne i}V'_{i,j} (d_{i,j} )\frac{(x_{i} -x_{j} )}{d_{i,j} }  $, $\sum _{j\ne i}p_{i,j} V'_{i,j} (d_{i,j} )$ \linebreak$\frac{(y_{i} -y_{j} )}{d_{i,j} }  $ contain at most $m_{i} $ non-zero terms, namely the terms with $d_{i,j} \le \lambda $. Definition \eqref{GrindEQ__5_17_} in conjunction with \eqref{GrindEQ__3_19_}, \eqref{GrindEQ__5_20_} and the fact that $p_{i,j} \ge 1$ for all $i,j=1,...,n$, $j\ne i$, implies the following estimate for all $w\in \Omega $ and $i=1,...,n$:
 
\begin{equation} \label{GrindEQ__5_21_} 
\begin{aligned}
&{\max \left(\left|\sum _{j\ne i}V'_{i,j} (d_{i,j} )\frac{(x_{i} -x_{j} )}{d_{i,j} }  \right|,\left|\sum _{j\ne i}p_{i,j} V'_{i,j} (d_{i,j} )\frac{(y_{i} -y_{j} )}{d_{i,j} }  \right|\right)} \\
&{\le \max \left(\sum _{j\ne i}\left|V'_{i,j} (d_{i,j} )\right|\left|\frac{x_{i} -x_{j} }{d_{i,j} } \right| ,\sum _{j\ne i}p_{i,j} \left|V'_{i,j} (d_{i,j} )\right|\left|\frac{y_{i} -y_{j} }{d_{i,j} } \right| \right)} \\ 
&\le \max \left(\sum _{j\ne i}\left|V'_{i,j} (d_{i,j} )\right| ,\sum _{j\ne i}\sqrt{p_{i,j} } \left|V'_{i,j} (d_{i,j} )\right| \right)\\
&\le \sum _{j\ne i}\sqrt{p_{i,j} } B_{i,j} \left(d_{i,j} \right) \le \sum _{j\ne i}\sqrt{p_{i,j} } B_{i,j} \left(\rho _{i,j} (H_{R} (w))\right) \\
&\le m_{i} \mathop{\max }\limits_{j\ne i} \left(\sqrt{p_{i,j} } B_{i,j} \left(\rho _{i,j} (H_{R} (w))\right)\right)
\end{aligned} 
\end{equation} 
Moreover, due to \eqref{GrindEQ__3_6_}, \eqref{GrindEQ__3_19_} and \eqref{GrindEQ__5_18_}, it holds that
\begin{equation} \label{GrindEQ__5_22_} 
\sum _{j\ne i}\kappa _{i,j} (d_{i,j} ) \le m_{i} \mathop{\max }\limits_{j\ne i} \left(c_{i,j} (\rho _{i,j} (H_{R} (w)))\right).   
\end{equation} 
Notice also that due to \eqref{GrindEQ__2_4_} and the facts that $v_{i} \in (0,v_{\max } )$ and $\cos (\theta )>\cos (\varphi )$, $\theta \in (-\varphi ,\varphi )$, $\varphi \in \left(-\frac{\pi }{2} ,\frac{\pi }{2} \right)$, we have that $v_{\max } v_{i} \cos (\theta _{i} )+v^{*} v_{\max } -2v^{*} v_{i} \ge v^{*} (v_{\max } -v_{i} )>0$, $i=1,...,n$. The previous inequality together with \eqref{GrindEQ__3_25_} implies that 
\begin{equation}   \label{GrindEQ__5_23_}
\frac{1}{q(v_{i} (t),\theta _{i} (t))} \le \frac{2(v_{\max } -v_{i} (t))v_{i}^{2} (t)}{v^{*} } \le \frac{2v_{\max }^{3} }{v^{*} } ,\,\,\textrm{for } t\ge0
\end{equation}

\noindent Using \eqref{GrindEQ__3_27_}, \eqref{GrindEQ__5_13_} and the facts that $\ell _{2} :\mathbb{R} _{+} \to [v^{*} ,v_{\max } )$, is non-decreasing and $\ell _{1} :\mathbb{R} _{+} \to (0,v^{*} ]$ is non-increasing, we obtain the following estimate,
\begin{equation} \label{GrindEQ__5_24_} 
a(v_{i} (t),\theta _{i} (t))\le \frac{bv_{\max } }{2\left(v_{\max } -\ell _{2} (H_{R} (w_{0} ))\right)^{2} \ell _{1} (H_{R} (w_{0} ))} .   
\end{equation} 
Moreover, due to the facts that $v_{i} \in \left(0,v_{\max } \right)$, $\theta _{i} \in \left(-\frac{\pi }{2} ,\frac{\pi }{2} \right)$ for $i=1,...,n$ and continuity of $f_{k} $, $k=1,2$, we have that 
\begin{align}
&\begin{aligned}
\left|f_{1} (v_{i} (t)\cos (\theta _{i} (t))-v^{*} )\right|\le \xi _{1} :=\mathop{\max }\limits_{x\in [-v^{*} ,v_{\max } -v^{*} ]} \left(|f_{1} (x)|\right)&\\
\textrm{for } i=1,...,n \,\,\textrm{and } t\ge0&
\end{aligned}\label{GrindEQ__5_25_}\\
&\begin{aligned}
\left|f_{2} (v_{i} (t)\sin (\theta _{i} (t)))\right|\le \xi _{2} :=\mathop{\max }\limits_{x\in [-v_{\max } ,v_{\max } ]} \left(|f_{2} (x)|\right)&\\
\textrm{for } i=1,...,n \,\,\textrm{and } t\ge0&
\end{aligned}\label{GrindEQ__5_26_}
\end{align}

\noindent Then, \eqref{GrindEQ__3_17_}, \eqref{GrindEQ__3_22_}, \eqref{GrindEQ__3_24_}, \eqref{GrindEQ__5_13_}, \eqref{GrindEQ__5_21_}, \eqref{GrindEQ__5_22_}, \eqref{GrindEQ__5_23_}, and \eqref{GrindEQ__5_25_} we obtain the following estimate for $i=1,...,n$ and $t\ge 0$ 
\begin{equation} \label{GrindEQ__5_27_} 
\begin{aligned}
&|F_{i} (t)|\le \frac{2v_{\max }^{3} }{v^{*} } \left(\xi _{1}\right.\\&\left. +m_{i} \left(g_{1} (v_{\max } )-g_{1} (0)\right)\mathop{\max }\limits_{j\ne i} \left(c_{i,j} (\rho _{i,j} (H_{R} (w_{0} )))\right)\right)\\ 
&+\frac{2v_{\max }^{3} }{v^{*} } m_{i} \mathop{\max }\limits_{j\ne i} \left(\sqrt{p_{i,j} } B_{i,j} (\rho _{i,j} (H_{R} (w_{0} )))\right)
\end{aligned} 
\end{equation} 
Estimate \eqref{GrindEQ__3_30_} with 
\[\begin{aligned} 
&Q_{1} (s):=\frac{2v_{\max }^{3} }{v^{*} } \left(\xi _{1} +\left(g_{1} (v_{\max } )-g_{1} (0)\right)\right.\\&\left.\times\mathop{\max }\limits_{i=1,...,n} \left(m_{i} \mathop{\max }\limits_{j\ne i} \left(c_{i,j} (\rho _{i,j} (s))\right)\right)\right) \\
&+\frac{2v_{\max }^{3} }{v^{*} } \mathop{\max }\limits_{i=1,...,n} \left(m_{i} \mathop{\max }\limits_{j\ne i} \left(\sqrt{p_{i,j} } B_{i,j} (\rho _{i,j} (s))\right)\right) \end{aligned}\] 
follows directly from inequality \eqref{GrindEQ__5_27_}. 

Using also \eqref{GrindEQ__3_26_} and the definition of $b>1-\frac{v^{*} }{v_{\max } } $, we also have that $\frac{1}{\beta (v_{i} (t),\theta _{i} (t))} \le \frac{(1-\cos (\varphi ))^{2} }{A} $ for $i=1,...,n$ and $t\ge 0$. The previous inequality, \eqref{GrindEQ__3_17_}, \eqref{GrindEQ__3_21_}, \eqref{GrindEQ__3_23_}, \eqref{GrindEQ__3_26_}, \eqref{GrindEQ__3_27_}, \eqref{GrindEQ__3_30_}, \eqref{GrindEQ__5_13_}, \eqref{GrindEQ__5_19_}, \eqref{GrindEQ__5_21_}, \eqref{GrindEQ__5_22_}, \eqref{GrindEQ__5_24_}, \eqref{GrindEQ__5_26_}, and the fact that $\ell _{2} (s)<v_{\max } $ for all $s\ge 0$,  give the following estimate for each $i=1,...,n$ and $t\ge 0$: 
\begin{equation} \label{GrindEQ__5_28_} 
\begin{aligned}
&|u_{i} (t)|\le \frac{v_{\max } (1-\cos (\varphi ))^{2} }{A} \\&\times\left(\xi _{2} +\frac{bv_{\max } Q_{1} (H_{R} (w_{0} ))}{2\left(v_{\max } -\ell _{2} (H_{R} (w_{0} ))\right)^{2} \ell _{1} (H_{R} (w_{0} ))}\right.\\& +\zeta _{i} \left(\eta _{i} (H_{R} (w_{0} ))\right)\Biggr)\\
&+\frac{v_{\max } (1-\cos (\varphi ))^{2} }{A} \frac{m_{i} }{v^{*} } \mathop{\max }\limits_{j\ne i} \left(\sqrt{p_{i,j} } B_{i,j} (\rho _{i,j} (H_{R} (w_{0} )))\right) \\ 
&+\frac{v_{\max } (1-\cos (\varphi ))^{2} }{A} \frac{m_{i} }{v^{*} } \left(g_{2} (v_{\max } )-g_{2} (-v_{\max } )\right)\\&\times\mathop{\max }\limits_{j\ne i} \left(c_{i,j} (\rho _{i,j} (H_{R} (w_{0} )))\right)
\end{aligned} 
\end{equation} 
Inequality \eqref{GrindEQ__3_31_} with
\[\begin{aligned}
&Q_{2} (s):=\frac{v_{\max } (1-\cos (\varphi ))^{2} }{A} \\&\times\left(\xi _{2} +\frac{bv_{\max } Q_{1} (s)}{2(v_{\max } -\ell _{2} (s))^{2} \ell _{1} (s)}  +\mathop{\max }\limits_{i=1,...,n} \left(\zeta _{i} (\eta _{i} (s))\right)\right) \\ 
&+\frac{v_{\max } (1-\cos (\varphi ))^{2} }{Av^{*} } \mathop{\max }\limits_{i=1,...,n} \left(m_{i} \mathop{\max }\limits_{j\ne i} \left(\sqrt{p_{i,j} } B_{i,j} (\rho _{i,j} (s))\right)\right) \\ 
&+\frac{v_{\max } (1-\cos (\varphi ))^{2} }{Av^{*} } \left(g_{2} (v_{\max } )-g_{2} (-v_{\max } )\right)\\&\times\mathop{\max }\limits_{i=1,...,n} \left(m_{i} \mathop{\max }\limits_{j\ne i} \left(c_{i,j} (\rho _{i,j} (s))\right)\right)
\end{aligned}\] 
is a direct consequence of \eqref{GrindEQ__5_28_}.

\noindent We show next that for all $i=1,\ldots,n$, $\frac{d}{dt}((v_{i} (t)\cos (\theta _{i} (t))-v^{*} )$ $f_{1} (v_{i} (t)\cos (\theta _{i} (t))-v^{*} ))$, and $\frac{d}{dt} (v_{i} (t)\sin (\theta _{i} (t))f_{2} (v_{i} (t)$ $\sin (\theta _{i} (t)))$ are bounded. Since $f_{k} \in C^{1} (\mathbb{R} )$, $k=1,2$ and due to the facts that $v_{i} \in(0,v_{\max })$, $\theta _{i} \in \left(-\frac{\pi }{2} ,\frac{\pi }{2} \right)$, we also have that
\begin{align}
&\begin{aligned}
\left|f_{1} ^{{'} } (v_{i} (t)\cos (\theta _{i} (t))-v^{*} )\right|\le \mathop{\max }\limits_{x\in [-v^{*} ,v_{\max } -v^{*} ]} \left(|f_{1} ^{{'} } (x)|\right)&\\\textrm{for } i=1,...,n \,\,\textrm{and } t\ge0&
\end{aligned}\label{GrindEQ__5_29_}\\
&\begin{aligned}
\left|f_{2} ^{{'} } (v_{i} (t)\sin (\theta _{i} (t)))\right|\le \mathop{\max }\limits_{x\in [-v_{\max } ,v_{\max } ]} \left(|f_{2} ^{{'} } (x)|\right)&\\
\textrm{for } i=1,...,n \,\,\textrm{and } t\ge0&
\end{aligned}\label{GrindEQ__5_30_}
\end{align}
 \nopagebreak
\noindent Estimates \eqref{GrindEQ__5_25_}, \eqref{GrindEQ__5_26_}, \eqref{GrindEQ__3_30_}, \eqref{GrindEQ__3_31_}, \eqref{GrindEQ__5_29_}, \eqref{GrindEQ__5_30_}, the facts that $v_{i} (t)\in (0,v_{\max } )$, $\theta _{i} (t)\in \left(-\varphi ,\varphi \right)$, for all $t\ge 0$, $i=1,...,n$,
imply that $\frac{d}{dt} \left((v_{i} (t)\cos (\theta _{i} (t))-v^{*} )f_{1} (v_{i} (t)\cos (\theta _{i} (t))-v^{*} )\right)$, $\frac{d}{dt} \left(v_{i} (t)\sin (\theta _{i} (t))f_{2} (v_{i} (t)\sin (\theta _{i} (t))\right)$ are bounded for all $i=1,...,n$. Moreover, by taking into account \eqref{GrindEQ__3_17_} and the fact that $v_{i} (t)\in (0,v_{\max } )$, $\theta _{i} (t)\in (-\varphi ,\varphi )$ for all $t\ge 0$, we have that for all $i=1,...,n$, and $t\ge 0$ 
\begin{align}
&\begin{aligned}
\left|g'_{1} \left(v_{i} (t)\cos (\theta _{i} (t))\right)\right|\le \mathop{\max }\limits_{x\in \left[0,v_{\max } \right]} \left(g'_{1} (x)\right)&\\\textrm{for } i=1,...,n \,\,\textrm{and } t\ge0&
\end{aligned}\label{GrindEQ__5_31_}\\
&\begin{aligned}
\left|g'_{2} \left(v_{i} (t)\sin (\theta _{i} (t))\right)\right|\le \mathop{\max }\limits_{x\in \left[-v_{\max } ,v_{\max } \right]} \left(g'_{2} (x)\right)&\\
\textrm{for } i=1,...,n \,\,\textrm{and } t\ge0&
\end{aligned}\label{GrindEQ__5_32_}
\end{align}

\noindent Inequalities \eqref{GrindEQ__5_31_}, \eqref{GrindEQ__5_32_} in conjunction with \eqref{GrindEQ__2_2_}, \eqref{GrindEQ__5_27_}, \eqref{GrindEQ__5_28_} and the fact that $v_{i} (t)\in (0,v_{\max } )$, $\theta _{i} (t)\in (-\varphi ,\varphi )$ for all $t\ge 0$, imply that $\frac{d}{dt} g_{1} (v_{i} (t)\cos (\theta _{i} (t)))$ and $\frac{d}{dt} g_{2} (v_{i} (t)$ $\sin (\theta _{i} (t)))$ are bounded for all $i=1,...,n$, $t\ge 0$. 

 Finally, we show that for  and for each  
\begin{equation}\label{GrindEQ__5_33_} 
\begin{aligned}
&\frac{d}{dt} \left(\sum _{j\ne i}\kappa _{i,j} \left(d_{i,j} (t)\right)\left(v_{j} (t)\cos (\theta _{j} (t))-v_{i} (t)\cos (\theta _{i} (t))\right)\right.\\
& \times\left(g_{1} \left(v_{j} (t)\cos (\theta _{j} (t))\right)-g_{1} \left(v_{i} (t)\cos (\theta _{i} (t))\right)\right) \Biggr)
\end{aligned}
\end{equation} 

\noindent is bounded. By virtue of \eqref{GrindEQ__5_27_}, \eqref{GrindEQ__5_28_}, and \eqref{GrindEQ__5_31_}, and the facts that $v_{i} (t)\in (0,v_{\max } )$, $\theta _{i} (t)\in (-\varphi ,\varphi )$, for all $t\ge 0$, and $i=1,...,n$, we have that for all $i,j=1,...,n$, with $j\neq i$, the terms $\frac{d}{dt} (v_{j} (t)\cos (\theta _{j} (t))-v_{i} (t)\cos (\theta _{i} (t)))$ and $\frac{d}{dt} (g_{1} (v_{j} (t)\cos (\theta _{j} (t)))-g_{1} (v_{i} (t)\cos (\theta _{i} (t))))$, are bounded. 

We finally show that $\frac{d}{dt} \kappa _{i,j} (d_{i,j} (t))=\kappa '_{i,j} (d_{i,j} (t))\dot{d}_{i,j} (t)$ is bounded for each $i,j=1,...,n$, $j\neq i$. We show first that $\dot{d}_{i,j} (t)$ is bounded for all $t\ge 0$, $i,j=1,...,n$ with $j\neq i$. Indeed, using the Cauchy-Schwarz inequality and definition \eqref{GrindEQ__2_5_}, the assumption that $p_{i,j} \ge 1$, for $i,j=1,...,n$, $j\neq i$, the facts that $v_{i} (t)\in (0,v_{\max } )$, $t\ge0$  and  $\sin(\theta)$ is increasing on the interval $\theta\in(-\varphi,\varphi)$, we get for $t\ge0$ and $i,j=1,...,n$, $j\neq i$:
\begin{equation} \label{GrindEQ__5_34_} 
\begin{aligned}
|\dot{d}_{i,j} (t)|&\le ( (v_{i} (t)\cos (\theta _{i} (t))-v_{j} (t)\cos (\theta _{j} (t)) )^{2} \\
&+p_{i,j} \left(v_{i} (t)\sin (\theta _{i} (t))-v_{j} (t)\sin (\theta _{j} (t))\right)^{2})^{\frac{1}{2}}   \\ 
&\le \sqrt{2} v_{\max } \sqrt{1+(2p_{i,j} -1)\sin ^{2} (\varphi )} =\delta _{i,j}  \end{aligned} 
\end{equation} 
Define the non-increasing functions
\begin{equation}\label{GrindEQ__5_35_}
\begin{aligned}
&\bar{\kappa }_{i,j} (s):=\max \left\{|\kappa '_{i,j} (d)|:s\le d\le \lambda \right\} \textrm{ for } s\in (L_{i,j},\lambda]\\
&i,j=1,...,n,\, j\ne i
\end{aligned}
\end{equation}

\noindent Then, from \eqref{GrindEQ__3_6_}, \eqref{GrindEQ__5_13_}, and \eqref{GrindEQ__5_35_} we get
\begin{equation} \label{GrindEQ__5_36_} 
\left|\kappa '(d_{i,j} (t))\right|\le \bar{\kappa }_{i,j} (\rho _{i,j} (H_{R} (w_{0} ))) , i,j=1,...,n,\, j\ne i.  
\end{equation} 
Inequalities \eqref{GrindEQ__5_34_} and \eqref{GrindEQ__5_36_} imply that $\frac{d}{dt} \kappa _{i,j} (d_{i,j} (t) =$ \linebreak$\kappa'_{i,j} (d_{i,j} (t))\dot{d}_{i,j} (t)$ is bounded for all $i,j=1,...,n,$ $j\ne i$. Combining, \eqref{GrindEQ__5_36_} and boundedness of $\frac{d}{dt} (v_{j} (t)\cos (\theta _{j} (t))-v_{i} (t)\cos (\theta _{i} (t)))$ and $\frac{d}{dt} (g_{1} (v_{j} (t)\cos (\theta _{j} (t)))-g_{1} (v_{i} (t)$\linebreak$\cos (\theta _{i} (t))))$, the fact that $v_i(t)\in(0,v_{\max})$, $\theta_i(t)\in(-\varphi,\varphi)$, for all $t\ge0$, and assumption \eqref{GrindEQ__3_17_}, we obtain that \eqref{GrindEQ__5_33_} is bounded. Finally, by virtue of \eqref{GrindEQ__5_27_}, \eqref{GrindEQ__5_28_}, and \eqref{GrindEQ__5_32_}, and the facts that $v_i(t)\in(0,v_{\max})$, $\theta_i(t)\in(-\varphi,\varphi)$, for all $t\ge0$, and $i=1,\ldots,n$, we also have that for all $i,j=1,...,n$, with $j\ne i$, the terms $\frac{d}{dt} (g_{2} (v_{j} (t)\sin (\theta _{j} (t)))$ $-g_{2} (v_{i} (t)\sin (\theta _{i} (t))))$, and $\frac{d}{dt} v_{j} (t)\sin (\theta _{j} (t))$ $-v_{i} (t)\sin (\theta _{i} (t)))$ are bounded. The latter, together with the fact that \linebreak$\frac{d}{dt} \kappa _{i,j} (d_{i,j} (t))$ is bounded for all $i,j=1,...,n$, $j\ne i$, and assumption \eqref{GrindEQ__3_17_}, implies that 
\begin{equation}\label{GrindEQ__5_37_} 
\begin{aligned}
&\frac{d}{dt} \left(\sum _{j\ne i}\kappa _{i,j} \left(d_{i,j} (t)\right)\left(v_{j} (t)\sin (\theta _{j} (t))-v_{i} (t)\sin (\theta _{i} (t))\right)\right.\\
&\left(g_{2} \left(v_{j} (t)\sin (\theta _{j} (t))\right)-g_{2} \left(v_{i} (t)\sin (\theta _{i} (t))\right)\right) \Biggr)
\end{aligned}
\end{equation}
is bounded for all $i=1,...,n$.

 Combining the facts that \eqref{GrindEQ__5_33_}, \eqref{GrindEQ__5_37_}, $\frac{d}{dt} ((v_{i} (t)\cos (\theta _{i} (t))-v^{*} )f_{1} (v_{i} (t)\cos (\theta _{i} (t))-v^{*} ))$,  and $\frac{d}{dt} (v_{i} (t)\sin (\theta _{i} (t))f_{2} (v_{i}(t) $\linebreak $ \sin (\theta _{i} (t)))$ are bounded for all $i=1,...,n$, we obtain that there exists $\bar{M}>0$  such that \eqref{GrindEQ__5_16_} holds.

Thus, from \eqref{GrindEQ__5_15_}, \eqref{GrindEQ__5_16_} and Barbalat's lemma (\cite{14}), we conclude that 

\noindent 
\begin{equation} \label{GrindEQ__5_38_} 
\mathop{\lim }\limits_{t\to +\infty } \left(\Delta (t)\right)=0.    
\end{equation} 
Since $\kappa _{i,j} (d)\ge 0$ for all $d>L_{i,j} $ and due to assumptions \eqref{GrindEQ__3_17_}, \eqref{GrindEQ__3_20_}, we have that the following inequalities hold for all $t\ge 0$, and $i=1,...,n$, 
\begin{align*}
&(v_i(t)\cos(\theta_i(t))-v^{*})f_1(v_i(t)\cos(\theta_i(t))-v^{*})\ge 0,\\
&(v_i(t)\sin(\theta_i(t)))f_2(v_i(t)\sin(\theta_i(t)))\ge 0,\\
&\sum _{j\ne i}\kappa _{i,j} \left(d_{i,j} (t)\right)\left(v_{j} (t)\sin (\theta _{j} (t))-v_{i} (t)\sin (\theta _{i} (t))\right)\\
&\left(g_{2} \left(v_{j} (t)\sin (\theta _{j} (t))\right)-g_{2} \left(v_{i} (t)\sin (\theta _{i} (t))\right)\right)\ge 0,\end{align*}\begin{align*}
&\sum _{j\ne i}\kappa _{i,j} \left(d_{i,j} (t)\right)\left(v_{j} (t)\cos (\theta _{j} (t))-v_{i} (t)\cos (\theta _{i} (t))\right)\\
& \times\left(g_{1} \left(v_{j} (t)\cos (\theta _{j} (t))\right)-g_{1} \left(v_{i} (t)\cos (\theta _{i} (t))\right)\right) \ge 0.
\end{align*}

\noindent Thus, we get from \eqref{GrindEQ__5_14_} that $0\le (v_{i} (t)\cos (\theta _{i} (t))-v^{*} )$ $f_{1} (v_{i} (t)\cos (\theta _{i} (t))-v^{*})\le \Delta (t)$, $0\le v_{i} (t)\sin (\theta _{i} (t))f_{2} (v_{i} (t)$ $\sin (\theta _{i} (t)))\le \Delta (t)$ for all $t\ge 0$, and $i=1,...,n$. The previous inequalities, in conjunction with \eqref{GrindEQ__5_13_}, \eqref{GrindEQ__5_38_} and assumption \eqref{GrindEQ__3_20_}, give \eqref{GrindEQ__3_28_} for the solution of the system \eqref{GrindEQ__2_2_}, \eqref{GrindEQ__3_21_}, \eqref{GrindEQ__3_22_}.

Finally, we show that \eqref{GrindEQ__3_29_} holds for the solutions  of \eqref{GrindEQ__2_2_}, \eqref{GrindEQ__3_21_}, \eqref{GrindEQ__3_22_} by exploiting Lemma 3 with $g(t)=v_i(t)$ and for $g(t)=\theta_i(t)$. Since \eqref{GrindEQ__3_28_} holds, it suffices to show that  and  are bounded. 

Using \eqref{GrindEQ__2_4_}, \eqref{GrindEQ__3_25_}, \eqref{GrindEQ__5_13_}, \eqref{GrindEQ__5_27_}, \eqref{GrindEQ__5_28_}, the fact that $v_{i} (t)\in (0,v_{\max } )$ and $\theta _{i} (t)\in (-\varphi ,\varphi )$ for all $t\ge 0$, $i=1,...,n$, inequality $v_{\max } v_{i} (t)\cos (\theta _{i} (t))+v^{*} v_{\max } -2v^{*} v_{i} (t)\ge v^{*} (v_{\max } -\ell _{2} (H_{R} (w_{0} )))$ and formula
\[\begin{aligned}
&{\left(v_{\max } v^{*} +v_{\max } \cos (\theta _{i} (t))-2v^{*} v_{i} (t)\right)^{2} \frac{d}{dt} \left(\frac{1}{q(v_{i} (t),\theta _{i} (t))} \right)=} \\ 
&=2\left(v_{\max } \sin (\theta _{i} (t))(v_{\max } -v_{i} (t))^{2} v_{i}^{3} (t)u_{i} (t)\right)\\&+4v^{*} v_{\max }^{2} (v_{\max } -v_{i} (t))v_{i} (t)-12v^{*} (v_{\max } -v_{i} (t))^{2} v_{i}^{2} (t)\\ 
&+2v_{\max } \cos (\theta _{i} (t))(v_{\max } -v_{i} (t))^{2} v_{i}^{2} (t)F_{i} (t)\\
&-4v_{\max } \cos (\theta _{i} (t))(v_{\max } -v_{i} (t))v_{i}^{3} (t)F_{i} (t) \end{aligned}\] 
we obtain that $\frac{d}{dt} \left(\frac{1}{q(v_{i} (t),\theta _{i} (t))} \right)$ is bounded for all $i=1,...,n$. 

We next prove that $\frac{d}{d\, t} \left(\sum _{j\ne i}V'_{i,j} (d_{i,j} (t))\frac{(x_{i} (t)-x_{j} (t))}{d_{i,j} (t)}  \right)$ is bounded for all $i=1,...,n$. Since \eqref{GrindEQ__3_2_} and \eqref{GrindEQ__5_13_} hold, it follows that $V'_{i,j} (d_{i,j} (t)),V''_{i,j} (d_{i,j} (t))$ are bounded for all $i,j=1,...,n$ with $j\ne i$. Moreover, due to \eqref{GrindEQ__5_13_}, \eqref{GrindEQ__5_20_}, \eqref{GrindEQ__5_21_}, \eqref{GrindEQ__5_34_}, it follows that
\begin{equation} \label{GrindEQ__5_39_} 
\begin{aligned}
&\frac{d}{d\, t} \left(\sum _{j\ne i}V'_{i,j} (d_{i,j} (t))\frac{(x_{i} (t)-x_{j} (t))}{d_{i,j} (t)}  \right)\\
&=\sum _{j\ne i}V''_{i,j} (d_{i,j} (t))\dot{d}_{i,j} (t)\frac{(x_{i} (t)-x_{j} (t))}{d_{i,j} (t)}  \\
&-\sum _{j\ne i}V'_{i,j} (d_{i,j} (t))\dot{d}_{i,j} (t)\frac{(x_{i} (t)-x_{j} (t))}{d_{i,j}^{2} (t)} \\
& +\sum _{j\ne i}V'_{i,j} (d_{i,j} (t))\frac{(v_{i} (t)\cos (\theta _{i} (t))-v_{j} (t)\cos (\theta _{j} (t)))}{d_{i,j} (t)}   \end{aligned} 
\end{equation} 
is bounded. Similarly, we prove that $\frac{d}{d\, t} (\sum _{j\ne i}p_{i,j} V'_{i,j} (d_{i,j} (t))$ $\frac{(y_{i} (t)-y_{j} (t))}{d_{i,j} (t)}  )$ is bounded for all $i=1,...,n$.  

Using \eqref{GrindEQ__3_6_}, \eqref{GrindEQ__3_17_}, \eqref{GrindEQ__3_30_}, \eqref{GrindEQ__3_31_}, \eqref{GrindEQ__5_22_}, \eqref{GrindEQ__5_31_}, \eqref{GrindEQ__5_32_}, \eqref{GrindEQ__5_34_} and boundedness of  (recall (5.34) and \eqref{GrindEQ__5_36_}),  $\frac{d}{dt} f_{1} (v_{i} (t)\cos (\theta _{i} (t))$ $-v^{*} )$, $\frac{d}{dt} f_{2} (v_{i} (t)\sin (\theta _{i} (t)))$, $\frac{d}{dt} g_{1} (v_{i} (t)\cos (\theta _{i} (t)))$, and $\frac{d}{dt} g_{2} $ $ (v_{i} (t)\sin (\theta _{i} (t)))$, we further obtain that $\frac{d}{dt} R_{i} (w(t))$ and $\frac{d}{dt} G_{i} (w(t))$, are also bounded for all $i=1,\ldots,n$, where $R_{i} (w(t))$ and $G_{i} (w(t))$, are defined in \eqref{GrindEQ__3_23_}, and \eqref{GrindEQ__3_24_}, respectively. Boundedness of $\frac{d}{dt} R_{i} (w(t))$ and $\frac{d}{dt} \left(\frac{1}{q(v_{i} (t),\theta _{i} (t))} \right)$, $i=1,\ldots,n$, in conjunction with \eqref{GrindEQ__3_22_}, \eqref{GrindEQ__5_13_}, \eqref{GrindEQ__5_21_}, \eqref{GrindEQ__5_22_}, \eqref{GrindEQ__5_23_}, \eqref{GrindEQ__5_25_}, \eqref{GrindEQ__5_34_}, \eqref{GrindEQ__5_39_}, imply that 
\[\begin{aligned}
&\dot{F}_{i} (t)=\frac{d}{dt} \left(\frac{1}{q(v_{i} (t),\theta _{i} (t))} \right)\\&\times\left(R_{i} (w(t))-\sum _{j\ne i}V'_{i,j} (d_{i,j} (t))\frac{x_{i} (t)-x_{j} (t)}{d_{i,j} (t)}  \right) \\ 
&+\frac{1}{q(v_{i} (t),\theta _{i} (t))} \\
&\times\left(\frac{d}{dt} R_{i} (w(t))-\frac{d}{dt} \left(\sum _{j\ne i}V'_{i,j} (d_{i,j} (t))\frac{x_{i} (t)-x_{j} (t)}{d_{i,j} (t)}  \right)\right) \end{aligned}\] 
is bounded for all $i=1,\ldots,n$.

Finally, since \eqref{GrindEQ__5_13_} holds, it follows that $U'(y_{i} (t))$ and $U''(y_{i} (t))$ are bounded for all $i=1,...,n$ which in conjunction with \eqref{GrindEQ__5_27_}, \eqref{GrindEQ__5_28_}, $v_{i} (t)\in (0,v_{\max } )$, $\theta _{i} (t)\in \left(-\varphi ,\varphi \right)$ for all $t\ge 0$, imply that $\frac{d}{dt} U'(y_{i} (t))$ is bounded. Using the facts that $v_{i} (t)\in (0,v_{\max } )$, $\theta _{i} (t)\in \left(-\varphi ,\varphi \right)$ for all $t\ge 0$, \eqref{GrindEQ__3_26_}, \eqref{GrindEQ__3_27_}, \eqref{GrindEQ__5_13_}, \eqref{GrindEQ__5_27_}, and \eqref{GrindEQ__5_28_}, inequalities$\frac{1}{\beta (v_{i} (t),\theta _{i} (t))} \le \frac{(1-\cos (\varphi ))^{2} }{A} $, $\frac{1}{v_{i} (t)} \le \frac{1}{\ell _{1} (H_{R} (w_{0} ))} $ , $\frac{1}{v_{\max } -v_{i} (t)} \le \frac{1}{v_{\max } -\ell _{2} (H_{R} (w_{0} ))} $ and formulas
\[\begin{aligned}
 &{\frac{d}{dt} a(v_{i} (t),\theta _{i} (t))} {=}{\frac{bv_{\max } \cos (\theta _{i} (t))u_{i} (t)}{2(v_{\max } -v_{i} (t))^{2} v_{i} (t)}  } \\  
& \qquad{-\frac{bv_{\max } \sin (\theta _{i} (t))F_{i} (t)}{2(v_{\max } -v_{i} (t))^{2} v_{i}^{2} (t)}+\frac{bv_{\max } \sin (\theta _{i} (t))F_{i} (t)}{(v_{\max } -v_{i} (t))^{3} v_{i} (t)} } \\[1em]
&\beta ^{2} (v_{i} (t),\theta _{i} (t))\frac{d}{dt} \left(\frac{1}{\beta (v_{i} (t),\theta _{i} (t))} \right)  {=}\\
&\left(\frac{(b-1)v_{i} (t)\sin (\theta _{i} (t))}{v_{\max } -v_{i} (t)} -\frac{2A\sin (\theta _{i} (t))}{(\cos (\theta _{i} (t))-\cos (\varphi ))^{3} } \right)u_{i} (t) \\ 
&{+\left(\frac{v^{*} +(b-1)v_{\max } \cos (\theta _{i} (t))}{(v_{\max } -v_{i} (t))^{2} } \right)F_{i} (t)} \end{aligned}\] 
it follows that $\frac{d}{dt}(\frac{1}{\beta(v_i(t),\theta_i(t))})$, and $\frac{d}{dt} (a(v_{i} (t),\theta _{i} (t))$ are bounded for all $t\ge 0$, $i=1,...,n$ as well. Moreover, using \eqref{GrindEQ__3_6_}, \eqref{GrindEQ__3_17_}, \eqref{GrindEQ__3_21_}, \eqref{GrindEQ__5_13_}, \eqref{GrindEQ__5_21_}, \eqref{GrindEQ__5_22_}, \eqref{GrindEQ__5_24_}, \eqref{GrindEQ__5_26_}, \eqref{GrindEQ__5_27_}, \eqref{GrindEQ__5_34_}, and boundedness of $\frac{d}{d\, t} \left(\sum _{j\ne i}p_{i,j} V'_{i,j} (d_{i,j} (t))\frac{(y_{i} (t)-y_{j} (t))}{d_{i,j} (t)}  \right)$, $\frac{d}{dt} G_{i} (w(t))$, $\frac{d}{dt} \left(\frac{1}{\beta (v_{i} (t),\theta _{i} (t))} \right)$, $\frac{d}{dt} (a(v_{i} (t),\theta _{i} (t))$, $\frac{d}{dt} U'(y_{i} (t))$,  and $\dot{F}_{i} (t)$, $i=1,...,n$ , we get that 
\[\begin{aligned}
&\dot{u}_{i} (t)=\frac{d}{dt} \left(\frac{v_{i} (t)}{\beta (v_{i} (t),\theta _{i} (t))} \right)\Biggl(G_{i} (w(t))-U'_{i} (y_{i} (t))\\
&\left.-a(v_{i} (t),\theta _{i} (t))F_{i} (t)-\sum _{j\ne i}p_{i,j} V'_{i,j} (d_{i,j} (t))\frac{\left(y_{i} (t)-y_{j} (t)\right)}{d_{i,j} (t)}  \right) \\
&+\frac{v_{i} (t)}{\beta (v_{i} (t),\theta _{i} (t))} \frac{d}{dt} \Biggl(G_{i} (w(t))-U'_{i} (y_{i} (t))\\
&\left.-a(v_{i} (t),\theta _{i} (t))F_{i} (t)-\sum _{j\ne i}p_{i,j} V'_{i,j} (d_{i,j} (t))\frac{\left(y_{i} (t)-y_{j} (t)\right)}{d_{i,j} (t)}  \right) \end{aligned}\] 
is bounded for all $t\ge 0$, $i=1,...,n$. This completes the proof. $\triangleleft$

\noindent \textbf{Proof of Theorem 2: }We show first certain properties of $k_{i} (w)$ in \eqref{GrindEQ__3_15_} that will be used in the proof of Theorem 2. For every $w\in \Omega $ the following inequalities hold 
\begin{equation} \label{GrindEQ__5_40_} 
k_{i} (w)v^{*} \ge \Lambda _{i} (w)\ge -k_{i} (w)(v_{\max } \cos (\theta _{i} )-v^{*} ) ,  i=1,...,n 
\end{equation} 
Indeed, inequality \eqref{GrindEQ__5_40_} is equivalent to the following inequality for  
\begin{equation} \label{GrindEQ__5_41_} 
\begin{aligned}
k_{i} (w)&\ge \max \left\{\frac{\Lambda _{i} (w)}{v^{*} } ,-\frac{\Lambda _{i} (w)}{v_{\max } \cos (\theta _{i} )-v^{*} } \right\}\\
&=\frac{\Lambda _{i} (w)}{v^{*} } +\frac{v_{\max } \cos (\theta _{i} )}{v^{*} (v_{\max } \cos (\theta _{i} )-v^{*} )} \max \left(0,-\Lambda _{i} (w)\right) 
\end{aligned}
\end{equation} 
Due to \eqref{GrindEQ__3_16_}, we have that for all $w\in \Omega $ , $\max \left(0,-\Lambda _{i} (w)\right)\le r\left(-\Lambda _{i} (w)\right)$ for all  $i=1,...,n$. The previous inequality, in conjunction with \eqref{GrindEQ__5_41_} and definition \eqref{GrindEQ__3_15_} implies inequality \eqref{GrindEQ__5_40_}.

 Moreover, for every $w\in \Omega $ it holds that
\begin{equation} \label{GrindEQ__5_42_} 
\mu _{2} \le k_{i} \left(w\right)\le \Theta _{i} \left(H(w)\right) 
\end{equation}

\noindent where $\Theta _{i} :\mathbb{R} _{+} \to \mathbb{R} _{+} $ for $i=1,...,n$ are the non-decreasing functions defined for $s\ge 0$ by the following formulas: 
\begin{equation} \label{GrindEQ__5_43_} 
\begin{aligned}
&\Theta _{i} (s):=\mu _{2} +\frac{m_{i} }{v^{*} } \left(\left(g_{1} (v_{\max } )-g_{1} (0)\right)\mathop{\max }\limits_{j\ne i} \left(c_{i,j} \left(\rho _{i,j} (s)\right)\right)\right.\\
&\left.+\mathop{\max }\limits_{j\ne i} \left(B_{i,j} \left(\rho _{i,j} (s)\right)\right)\right) \\ 
&+\left(v_{\max } \left(1-\cos (\varphi )\right)\left(A+s\cos (\varphi )\right)\right)\\
&\times\max \left(r(z):\left|z\right|\le m_{i} \left(\mathop{\max }\limits_{j\ne i} \left(\sqrt{p_{i,j} } B_{i,j} (\rho _{i,j} (s))\right)\right.\right.\\
&\left.+\left(g_{1} (v_{\max } )-g_{1} (0)\right)\mathop{\max }\limits_{j\ne i} \left(c_{i,j} (\rho _{i,j} (s))\right)\right) \\
&\times  \left(Av^{*} \left(v_{\max } -v^{*} \right)+v^{*} \left(v_{\max } \cos (\varphi )-v^{*} \right)\left(1-\cos (\varphi )\right)s\right)^{-1}
 \end{aligned} 
\end{equation} 
\noindent where \textit{$\rho _{i,j} :\mathbb{R} _{+} \to (L_{i,j} ,\lambda ]$} for $i,j=1,...,n$, $j\ne i$ are the non-increasing functions involved in \eqref{GrindEQ__3_10_}, $H(w)$ is given by \eqref{GrindEQ__3_9_}, $r(x)$ satisfies \eqref{GrindEQ__3_16_}, and $B_{i,j} (s)$, $c_{i,j} (s)$, for $i,j=1,...,n$, $j\ne i$ are defined in \eqref{GrindEQ__5_17_} and \eqref{GrindEQ__5_18_}, respectively. In order to show \eqref{GrindEQ__5_42_} notice first that due to definition of $m_{i} $, \eqref{GrindEQ__5_17_}, \eqref{GrindEQ__3_10_}, and \eqref{GrindEQ__5_20_}, we have (with similar arguments that led to the derivation of \eqref{GrindEQ__5_21_}) that
\begin{equation} \label{GrindEQ__5_44_} 
\begin{aligned}& {\max \left(\left|\sum _{j\ne i}V'_{i,j} (d_{i,j} )\frac{(x_{i} -x_{j} )}{d_{i,j} }  \right|,\left|\sum _{j\ne i}p_{i,j} V'_{i,j} (d_{i,j} )\frac{(y_{i} -y_{j} )}{d_{i,j} }  \right|\right)} \\& \le \sum _{j\ne i}\sqrt{p_{i,j} } B_{i,j} \left(\rho _{i,j} (H(w))\right)\\ &\le m_{i} \mathop{\max }\limits_{j\ne i} \left(\sqrt{p_{i,j} } B_{i,j} \left(\rho _{i,j} (H(w))\right)\right)\end{aligned} 
\end{equation}

\noindent Then, from definitions \eqref{GrindEQ__3_14_}, \eqref{GrindEQ__5_17_}, \eqref{GrindEQ__5_18_} and \eqref{GrindEQ__3_6_}, \eqref{GrindEQ__5_44_}, and the facts that $g_{1} (\cdot )$ satisfies \eqref{GrindEQ__3_17_}, $v_{i} \in (0,v_{\max } )$ and $\theta _{i} \in \left(-\varphi ,\varphi \right)$, we have that
\begin{equation} \label{GrindEQ__5_45_} 
\begin{aligned}\left|\Lambda _{i} (w)\right|\le& m_{i} \mathop{\max }\limits_{j\ne i} \left(\sqrt{p_{i,j} } B_{i,j} \left(\rho _{i,j} (H(w))\right)\right) \\  &+m_{i} \mathop{\max }\limits_{j\ne i} \left(c_{i,j} \left(\rho _{i,j} (H(w))\right)\right)\left(g_{1} (v_{\max } )-g_{1} (0)\right)\\&\qquad \qquad \qquad i=1,...,n.    \end{aligned}  
\end{equation} 
In addition, definition \eqref{GrindEQ__3_9_} implies that the following estimate holds for all  and  
\begin{equation} \label{GrindEQ__5_46_} 
\begin{aligned}
&\frac{1}{v_{\max } \cos \left(\theta _{i} \right)-v^{*} } \le\\
&\frac{A+\left(1-\cos (\varphi )\right)H(w)}{A\left(v_{\max } -v^{*} \right)+\left(v_{\max } \cos (\varphi )-v^{*} \right)\left(1-\cos (\varphi )\right)H(w)} .   
\end{aligned}
\end{equation} 
Inequality, \eqref{GrindEQ__5_42_} with $\Theta _{i} :\mathbb{R} _{+} \to \mathbb{R} _{+} $, $i=1,...,n$, defined by \eqref{GrindEQ__5_43_} is a direct consequence of  \eqref{GrindEQ__5_41_}, \eqref{GrindEQ__5_44_}, \eqref{GrindEQ__5_45_}, \eqref{GrindEQ__5_46_} and definition \eqref{GrindEQ__3_15_}.

Let $w_{0} \in \Omega $ and consider the unique solution $w(t)$ of the initial value problem \eqref{GrindEQ__2_2_}, \eqref{GrindEQ__3_11_}, \eqref{GrindEQ__3_12_} with initial condition $w(0)=w_0$. Using the fact that the set $\Omega $ is open (recall definitions \eqref{GrindEQ__2_3_}, \eqref{GrindEQ__2_7_}), we conclude that there exists $t_{\max } \in (0,+\infty ]$ such that the solution $w(t)$ of \eqref{GrindEQ__2_2_}, \eqref{GrindEQ__3_11_}, \eqref{GrindEQ__3_12_} is defined on $\left[0,t_{\max } \right)$ and satisfies $w(t)\in \Omega $ for all $t\in \left[0,t_{\max } \right)$. Furthermore, if $t_{\max } <+\infty $ then there exists an increasing sequence of times $\left\{\, t_{i} \in \left[0,t_{\max } \right)\, :\, i=1,2,...\right\}$ with $\mathop{\lim }\limits_{i\to +\infty } \left(t_{i} \right)=t_{\max } $ and either $\mathop{\lim }\limits_{i\to +\infty } \left(dist\left(w(t_{i} ),\partial \Omega \right)\right)=0$ or $\mathop{\lim }\limits_{i\to +\infty } \left(\left|w(t_{i} )\right|\right)=+\infty $.

We show that the solution $w(t)$ of the initial value problem \eqref{GrindEQ__2_2_}, \eqref{GrindEQ__3_11_}, \eqref{GrindEQ__3_12_} satisfies $w(t)\in \Omega $ for all $t\ge 0$. Using \eqref{GrindEQ__2_2_}, definition \eqref{GrindEQ__3_9_} and the fact that $\sum _{i=1}^{n}\sum _{j\ne i}^{}V'_{i,j} (d_{i,j} )$ $\frac{(x_{i} -x_{j} )}{d_{i,j} }   =0$ we obtain that
\begin{equation} \label{GrindEQ__5_46b_} 
\begin{aligned}
&\nabla H(w)\dot{w}=\\
&\sum _{i=1}^{n}\left(v_{i} \cos (\theta _{i} )-v^{*} \right)\left(F_{i} \cos (\theta _{i} )+\sum _{j\ne i}V'_{i,j} (d_{i,j} )\frac{(x_{i} -x_{j} )}{d_{i,j} }  \right)  \\ 
&+\sum _{i=1}^{n}v_{i} \sin (\theta _{i} )\left(b\sin (\theta _{i} )F_{i} +\left(v^{*} +\frac{A}{v_{i} \left(\cos (\theta _{i} )-\cos (\varphi )\right)^{2} }\right.\right.\\
& +(b-1)v_{i} \cos (\theta _{i} )\Biggr)u_{i} \Biggr)  \\
&+\sum _{i=1}^{n}v_{i} \sin (\theta _{i} )\left(U'_{i} (y_{i} )+\sum _{j\ne i}p_{i,j} V'_{i,j} (d_{i,j} )\frac{(y_{i} -y_{j} )}{d_{i,j} }  \right)  \end{aligned} 
\end{equation} 
Finally, \eqref{GrindEQ__3_3_}, \eqref{GrindEQ__3_7_}, \eqref{GrindEQ__3_11_}, \eqref{GrindEQ__3_12_}, \eqref{GrindEQ__3_15_}, \eqref{GrindEQ__5_46b_} and the fact that $k_{i} (w)\ge \mu _{2} $ (recall \eqref{GrindEQ__5_42_}), imply that for all $i=1,...,n$, $w\in \Omega $, the following inequality holds 
\begin{equation} \label{GrindEQ__5_47_} 
\begin{aligned} 
&{\nabla H(w)\dot{w}}  {=}  {-\sum _{i=1}^{n}k_{i} (w)\left(v_{i} \cos (\theta _{i} )-v^{*} \right)^{2}  -\mu _{1} \sum _{i=1}^{n}v_{i}^{2} \sin ^{2} (\theta _{i} ) } \\ 
&-\frac{1}{2} \sum _{i=1}^{n}\sum _{j\ne i}\kappa _{i,j} \left(d_{i,j} (t)\right)\left(v_{j} (t)\sin (\theta _{j} (t))-v_{i} (t)\sin (\theta _{i} (t))\right)\\
&\times\left(g_{2} \left(v_{j} (t)\sin (\theta _{j} (t))\right)-g_{2} \left(v_{i} (t)\sin (\theta _{i} (t))\right)\right)   \\
& -\frac{1}{2} \sum _{i=1}^{n}\sum _{j\ne i}\kappa _{i,j} \left(d_{i,j} (t)\right)\left(v_{j} (t)\cos (\theta _{j} (t))-v_{i} (t)\cos (\theta _{i} (t))\right)\\
&\times\left(g_{1} \left(v_{j} (t)\cos (\theta _{j} (t))\right)-g_{1} \left(v_{i} (t)\cos (\theta _{i} (t))\right)\right)\\
&\leq{-\sum _{i=1}^{n}\mu_2\left(v_{i} \cos (\theta _{i} )-v^{*} \right)^{2}  -\mu _{1} \sum _{i=1}^{n}v_{i}^{2} \sin ^{2} (\theta _{i} ) } \\
&-\frac{1}{2} \sum _{i=1}^{n}\sum _{j\ne i}\kappa _{i,j} \left(d_{i,j} (t)\right)\left(v_{j} (t)\sin (\theta _{j} (t))-v_{i} (t)\sin (\theta _{i} (t))\right)\\
&\times\left(g_{2} \left(v_{j} (t)\sin (\theta _{j} (t))\right)-g_{2} \left(v_{i} (t)\sin (\theta _{i} (t))\right)\right)   \\
& -\frac{1}{2} \sum _{i=1}^{n}\sum _{j\ne i}\kappa _{i,j} \left(d_{i,j} (t)\right)\left(v_{j} (t)\cos (\theta _{j} (t))-v_{i} (t)\cos (\theta _{i} (t))\right)\\
&\times\left(g_{1} \left(v_{j} (t)\cos (\theta _{j} (t))\right)-g_{1} \left(v_{i} (t)\cos (\theta _{i} (t))\right)\right)\\
  \end{aligned} 
\end{equation}

\noindent It follows from \eqref{GrindEQ__3_17_} and \eqref{GrindEQ__5_47_} that for all $w\in \Omega $ 
\begin{equation} \label{GrindEQ__5_48_} 
\nabla H(w)\dot{w}\le 0.                                                              
\end{equation}

\noindent Since $w(t)\in \Omega $ for all $t\in \left[0,t_{\max } \right)$, it follows from \eqref{GrindEQ__2_3_}, \eqref{GrindEQ__2_7_} that $v_{i} (t)\in \left(0,v_{\max } \right)$, $\theta _{i} (t)\in \left(-\varphi ,\varphi \right)$ for all $t\in \left[0,t_{\max } \right)$ and $i=1,...,n$. Thus, \eqref{GrindEQ__2_2_} implies that $0\le \dot{x}_{i} (t)\le v_{\max } $ for all $t\in \left[0,t_{\max } \right)$ and $i=1,...,n$. Moreover, inequality \eqref{GrindEQ__5_48_} gives that 
\begin{equation} \label{GrindEQ__5_49_}
H(w(t))\le H(w_{0} ), \textrm{ for all }t\in \left[0,t_{\max } \right).         
\end{equation}

\noindent Consequently, we obtain from \eqref{GrindEQ__3_10_}, \eqref{GrindEQ__5_49_} for all $t\in \left[0,t_{\max } \right)$ and $i,j=1,...,n$, $j\ne i$:
\begin{equation} \label{GrindEQ__5_50_} 
\begin{array}{l} {\left|\theta _{i} (t)\right|\le \omega \left(H(w_{0} )\right)<\varphi ,\left|y_{i} (t)\right|\le \eta _{i} \left(H(w_{0} )\right)<a,} \\ {x_{i} (0)\le x_{i} (t)\le x_{i} (0)+v_{\max } t\; ,\; d_{i,j} (t)\ge \rho _{i,j} \left(H(w_{0} )\right)>L_{i,j} } \end{array} 
\end{equation}

\noindent From \eqref{GrindEQ__5_42_} and \eqref{GrindEQ__5_49_} we finally get that
\begin{equation}\label{GrindEQ__5_51_} 
\begin{aligned}
\mu _{2} <k_{i} \left(w(t)\right)\le M_{i} :=\Theta _{i} \left(H(w_{0} )\right),& \textrm{ for all }t\in \left[0,t_{\max } \right)\\& \textrm{ and  }i=1,...,n      
\end{aligned}
\end{equation}
 
\noindent where $\Theta _{i} :[0,+\infty )\to \mathbb{R} _{+} $ , $i=1,...,n$ are increasing functions defined by \eqref{GrindEQ__5_43_}.

Inequalities \eqref{GrindEQ__5_40_}, the facts that $v_{i} (t)\in \left(0,v_{\max } \right)$, $t\in \left[0,t_{\max } \right)$, $k_{i} (w)>0$, for all $w\in \Omega $, and definitions \eqref{GrindEQ__3_12_} and \eqref{GrindEQ__3_14_}, imply the following differential inequalities for all $t\in \left[0,t_{\max } \right)$ and $i=1,...,n$:
\begin{equation} \label{GrindEQ__5_52_} 
k_{i} \left(w(t)\right)\left(v_{\max } -v_{i} (t)\right)\ge \dot{v}_{i} (t)\ge -k_{i} \left(w(t)\right)v_{i} (t).     
\end{equation} 
Differential inequalities \eqref{GrindEQ__5_52_} in conjunction with inequalities \eqref{GrindEQ__5_51_} imply that the following estimates hold for all $t\in \left[0,t_{\max } \right)$ and $i=1,...,n$:
\begin{equation} \label{GrindEQ__5_53_} 
\begin{aligned}
&v_{i} (0)\exp \left(-M_{i} \, t\right)+\left(1-\exp \left(-M_{i} \, t\right)\right)v_{\max } \ge v_{i} (t)\\
&\ge v_{i} (0)\exp \left(-M_{i} \, t\right).   
\end{aligned}               
\end{equation} 
Suppose that $t_{\max } <+\infty $. Inequalities \eqref{GrindEQ__5_49_}, \eqref{GrindEQ__5_50_}, \eqref{GrindEQ__5_53_} and definitions \eqref{GrindEQ__2_3_}, \eqref{GrindEQ__2_7_} imply that for every increasing sequence of times $\left\{\, t_{i} \in \left[0,t_{\max } \right)\, :\, i=1,2,...\right\}$ with $\mathop{\lim }\limits_{i\to +\infty } \left(t_{i} \right)=t_{\max } $ we cannot have  or $\mathop{\lim }\limits_{i\to +\infty } \left(\left|w(t_{i} )\right|\right)=+\infty $. Thus, we must have $t_{\max}=+\infty$. 

Since $w(t)\in \Omega $ of all $t\ge 0$ for system \eqref{GrindEQ__2_2_}, \eqref{GrindEQ__3_11_}, \eqref{GrindEQ__3_12_}, and due to \eqref{GrindEQ__2_2_}, \eqref{GrindEQ__5_42_}, \eqref{GrindEQ__5_49_}, and \eqref{GrindEQ__5_52_}, we have that 
\begin{equation} \label{GrindEQ__5_54_} 
\left|\dot{v}_{i} (t)\right|=|F_{i} (t)|\le \Theta _{i} (H(w_{0} ))v_{\max } ,\,i=1,...,n
\end{equation}

\noindent which implies that estimate \eqref{GrindEQ__3_34_} holds with 
\[Q_{3} (s)=\mathop{\max }\limits_{i=1,...,n} \left(\Theta _{i} (s)\right)v_{\max } . \] 
Next, using the definition of $b>1-\frac{v^{*} }{v_{\max } } $ and the fact that $v_{i} \in (0,v_{\max } )$, $\cos (\theta _{i} )>\cos (\varphi )$, $\theta \in (-\varphi ,\varphi )$, we have the following estimate
\begin{equation} \label{GrindEQ__5_55_} 
\begin{aligned}
&v^{*} +\frac{A}{v_{i} \left(\cos (\theta _{i} )-\cos (\varphi )\right)^{2} } +v_{i} \cos (\theta _{i} )(b-1) \\ 
&>v^{*} \left(1-\frac{v_{i} \cos (\theta _{i} )}{v_{\max } } \right)+\frac{A}{v_{i} \left(\cos (\theta _{i} )-\cos (\varphi )\right)^{2} } \\
&\ge \frac{A}{v_{\max } \left(1-\cos (\varphi )\right)^{2} } 
\end{aligned} 
\end{equation}

\noindent From \eqref{GrindEQ__2_2_}, \eqref{GrindEQ__5_18_}, \eqref{GrindEQ__5_19_}, \eqref{GrindEQ__5_44_}, \eqref{GrindEQ__5_49_}, \eqref{GrindEQ__5_50_}, \eqref{GrindEQ__5_45_}, \eqref{GrindEQ__5_51_}, \eqref{GrindEQ__5_43_}, \eqref{GrindEQ__5_54_} and \eqref{GrindEQ__5_55_}, we obtain the following estimate for all $t\ge 0$ and $i=1,...,n$
\begin{equation} \label{GrindEQ__5_56_} 
\begin{aligned}
{\left|\dot{\theta }_{i} (t)\right|} & {=} \left|u_{i} (t)\right|\le \frac{v_{\max } (1-\cos (\varphi ))^{2} }{A} \left(\left(\mu _{1} +b\Theta _{i} (H(w_{0} ))\right)v_{\max }\right.\\&\left. +\zeta _{i} (\eta _{i} (H(w_{0} ))\right) \\  
& {+\frac{m_{i} v_{\max } (1-\cos (\varphi ))^{2} }{A} \mathop{\max }\limits_{j\ne i} \left(\sqrt{p_{i,j} } B_{i,j} (\rho _{i,j} (H(w_{0} )))\right)} \\ 
& +\frac{m_{i} v_{\max } (1-\cos (\varphi ))^{2} }{A} \left(g_{2} (v_{\max } )-g_{2} (-v_{\max } )\right)\\
&\times\mathop{\max }\limits_{j\ne i} \left(c_{i,j} (\rho _{i,j} (H(w_{0} )))\right)
\end{aligned} 
\end{equation} 
Inequality \eqref{GrindEQ__3_35_} with
\[\begin{aligned}
&Q_{4} (s):=\frac{v_{\max } (1-\cos (\varphi ))^{2} }{A} \left(\left(\mu _{1} +b\mathop{\max }\limits_{i=1,...,n} \left(\Theta _{i} (s)\right)\right)v_{\max } \right.\\
&\left.+\mathop{\max }\limits_{i=1,...,n} \left(\zeta _{i} (\eta _{i} (s))\right)\right) \\ 
&+\frac{v_{\max } (1-\cos (\varphi ))^{2} }{Av^{*} } \mathop{\max }\limits_{i=1,...,n} \left(m_{i} \mathop{\max }\limits_{j\ne i} \left(\sqrt{p_{i,j} } B_{i,j} (\rho _{i,j} (s))\right)\right)\\ 
&+\frac{v_{\max } (1-\cos (\varphi ))^{2} }{Av^{*} } \left(g_{2} (v_{\max } )-g_{2} (-v_{\max } )\right)\\&\times\mathop{\max }\limits_{i=1,...,n} \left(m_{i} \mathop{\max }\limits_{j\ne i} \left(c_{i,j} (\rho _{i,j} (s))\right)\right) \end{aligned}\]

\noindent is a direct consequence of \eqref{GrindEQ__5_56_}.

We show next that \eqref{GrindEQ__3_32_} holds for system \eqref{GrindEQ__2_2_}, \eqref{GrindEQ__3_11_}, \eqref{GrindEQ__3_12_}. Define
\begin{equation}\label{GrindEQ__5_57_}
\begin{aligned}
& \Gamma (t):=\\
&\sum _{i=1}^{n}\mu _{2} \left(v_{i} (t)\cos (\theta _{i} (t))-v^{*} \right)^{2}  +\sum _{i=1}^{n}\mu _{1} \left(v_{i} (t)\sin (\theta _{i} (t))\right)^{2}  \\ 
& +\frac{1}{2} \sum _{i=1}^{n}\sum _{j\ne i}\kappa _{i,j} \left(d_{i,j} (t)\right)\left(v_{j} (t)\sin (\theta _{j} (t))-v_{i} (t)\sin (\theta _{i} (t))\right)\\&\times\left(g_{2} \left(v_{j} (t)\sin (\theta _{j} (t))\right)-g_{2} \left(v_{i} (t)\sin (\theta _{i} (t))\right)\right)   \\ 
&+\frac{1}{2} \sum _{i=1}^{n}\sum _{j\ne i}\kappa _{i,j} \left(d_{i,j} (t)\right)\left(v_{j} (t)\cos (\theta _{j} (t))-v_{i} (t)\cos (\theta _{i} (t))\right)\\&\times\left(g_{1} \left(v_{j} (t)\cos (\theta _{j} (t))\right)-g_{1} \left(v_{i} (t)\cos (\theta _{i} (t))\right)\right)  
\end{aligned}
\end{equation} 

\noindent Definitions \eqref{GrindEQ__5_57_}, \eqref{GrindEQ__5_47_}, and \eqref{GrindEQ__5_48_} imply that $\Gamma (t)=$\linebreak $-\frac{d}{dt} H(w(t))\ge 0$  for all $t\ge 0$. Therefore, since $H(w)\ge 0$ for all $w\in \Omega $, we obtain:
\begin{equation} \label{GrindEQ__5_58_} 
\int _{0}^{\infty }\Gamma (t)dt \le H(w_{0} ).                                                               
\end{equation}

\noindent We will show next that there exists $\tilde{M}>0$ such that 
\begin{equation}\label{GrindEQ__5_59_}
\left|\frac{d}{dt} \left(\Gamma (t)\right)\right|\le \tilde{M},\textrm{ for all  }t\ge 0
\end{equation} 

\noindent Notice that, due to \eqref{GrindEQ__5_54_}, \eqref{GrindEQ__5_56_} and the fact that $v_{i} (t)\in (0,v_{\max } )$, $\theta _{i} (t)\in (-\varphi ,\varphi )$ for all $t\ge 0$, and , we have that $\frac{d}{dt} \left(v_{i} (t)\cos (\theta _{i} (t))-v^{*} \right)^{2} $ and $\frac{d}{dt} \left(v_{i} (t)\sin (\theta _{i} (t))\right)^{2} $ are bounded. Moreover, from \eqref{GrindEQ__5_18_}, \eqref{GrindEQ__5_35_} and \eqref{GrindEQ__5_50_} we have that
\begin{equation} \label{GrindEQ__5_60_} 
\left|\kappa '_{i,j} (d_{i,j} (t))\right|\le \bar{\kappa }_{i,j} (\rho _{i,j} (H(w_{0} ))), i,j=1,...,n,\, \, \, j\ne i 
\end{equation} 
\begin{equation} \label{GrindEQ__5_61_} 
0\le \kappa _{i,j} (d_{i,j} (t))\le c_{i,j} (\rho _{i,j} (H(w_{0} ))).   
\end{equation}

\noindent Inequalities \eqref{GrindEQ__5_34_}, \eqref{GrindEQ__5_60_}, \eqref{GrindEQ__5_61_} imply that $\frac{d}{dt} \kappa _{i,j} (d_{i,j} (t))$ is bounded for all $i,j=1,...,n$, $j\ne i$. From the latter, together with boundedness of $\frac{d}{dt} g_{1} (v_{i} (t)\cos (\theta _{i} (t))$ and $\frac{d}{dt} g_{2} (v_{i} (t)\sin (\theta _{i} (t))$ (recall \eqref{GrindEQ__5_31_}, \eqref{GrindEQ__5_32_}), the fact that $v_{i} (t)\in (0,v_{\max } )$, $\theta _{i} (t)\in (-\varphi ,\varphi )$ for all $t\ge 0$, and , we obtain that there exists $\tilde{M}>0$ such that \eqref{GrindEQ__5_59_} holds. 
 
Thus, from Lemma 3, \eqref{GrindEQ__5_58_} and \eqref{GrindEQ__5_59_} we conclude that $\mathop{\lim }\limits_{t\to +\infty } \left(\Gamma (t)\right)=0$. Since $\kappa _{i,j} (d)\ge 0$ for all $d>L_{i,j} $, $\Gamma (t)\ge 0$ for all $t\ge 0$, and due to assumption \eqref{GrindEQ__3_20_}, we have that \eqref{GrindEQ__3_32_} holds for the solution $w(t)$ of \eqref{GrindEQ__2_2_}, \eqref{GrindEQ__3_11_}, \eqref{GrindEQ__3_12_}.

Finally, we show that \eqref{GrindEQ__3_33_} holds for the solution $w(t)$ of \eqref{GrindEQ__2_2_}, \eqref{GrindEQ__3_11_}, \eqref{GrindEQ__3_12_} by exploiting Lemma 3 with $g(t)=v_{i} (t)$ and $g(t)=\theta _{i} (t)$ for $i=1,...,n$. Since \eqref{GrindEQ__3_32_} holds, it suffices to show that $\dot{u}_i(t)$ and $\dot{F}_i(t)$ are bounded. 

Boundedness of $\frac{d}{d\, t} \left(\sum _{j\ne i}V'_{i,j} (d_{i,j} (t))\frac{(x_{i} (t)-x_{j} (t))}{d_{i,j} (t)}  \right)$ and \linebreak $\frac{d}{d\, t} \left(\sum _{j\ne i}p_{i,j} V'_{i,j} (d_{i,j} (t))\frac{(y_{i} (t)-y_{j} (t))}{d_{i,j} (t)}  \right)$ for each  follows directly from \eqref{GrindEQ__5_34_}, \eqref{GrindEQ__5_44_}, \eqref{GrindEQ__5_50_}, and formula \eqref{GrindEQ__5_39_}.  Since \eqref{GrindEQ__3_2_} and \eqref{GrindEQ__5_50_} hold, it follows that $V'_{i,j} (d_{i,j} (t))$, $V''_{i,j} (d_{i,j} (t))$ are bounded for all $i,j=1,...,n$ with $j\ne i$. The previous bounds together with \eqref{GrindEQ__5_34_}, \eqref{GrindEQ__5_31_}, \eqref{GrindEQ__5_32_}, \eqref{GrindEQ__5_60_}, and \eqref{GrindEQ__5_61_} imply that $\frac{d}{dt} \left(\Lambda _{i} (w(t))\right)$ for  are also bounded. 

We show next that $\frac{d}{dt} k_{i} (w(t))$ is bounded. We have
\begin{equation} \label{GrindEQ__5_62_} 
\begin{aligned}
\frac{d}{d\, t}& k_{i} (w(t))  {=}  
\frac{1}{v^{*} } \left(1-r'\left(-\Lambda _{i} (w(t))\right)\right)\frac{d}{d\, t} \left(\Lambda _{i} (w(t))\right)  \\
  & {-\frac{1}{v_{\max } \cos (\theta _{i} (t))-v^{*} } r'\left(-\Lambda _{i} (w(t))\right)\frac{d}{d\, t} \left(\Lambda _{i} (w(t))\right)} \\ 
 & {+\frac{v_{\max } \sin (\theta _{i} (t))u_{i} (t)}{\left(v_{\max } \cos (\theta _{i} (t))-v^{*} \right)^{2} } r\left(-\Lambda _{i} (w(t))\right)} \end{aligned} 
\end{equation} 
Boundedness of $\frac{d}{dt} \left(\Lambda _{i} (w(t))\right)$,  \eqref{GrindEQ__5_45_}, \eqref{GrindEQ__5_49_}, \eqref{GrindEQ__5_56_}, \eqref{GrindEQ__5_62_}, and inequalities
\[r\left(-\Lambda _{i} (w(t))\right)\le \max \left(r(z):\left|z\right|\le \psi _{i} \right),   \] 
\[\left|r'\left(-\Lambda _{i} (w(t))\right)\right|\le \max \left(\left|r'(z)\right|:\left|z\right|\le \psi _{i} \right),    \] 
where
\[
\begin{aligned}
\psi _{i} &=m_{i} \left(\left(g_{1} (v_{\max } )-g_{1} (0)\right)\mathop{\max }\limits_{j\ne i} \left(c_{i,j} (\rho _{i,j} (H(w_{0} )))\right)\right.\\
&\left.+\mathop{\max }\limits_{j\ne i} \left(\sqrt{p_{i,j} } B_{i,j} (\rho _{i,j} (H(w_{0} )))\right)\right) ,
\end{aligned}
\]

\noindent imply that each $\frac{d}{dt} \left(k_{i} (w(t))\right)$  is bounded for $i=1,...,n$.

Finally, from boundedness of $\frac{d}{dt} \left(\Lambda _{i} (w(t))\right)$, $\frac{d}{dt} \left(k_{i} (w(t))\right)$, \eqref{GrindEQ__5_42_}, \eqref{GrindEQ__5_45_}, \eqref{GrindEQ__5_49_}, \eqref{GrindEQ__5_54_}, \eqref{GrindEQ__5_56_}, inequality $\frac{1}{\cos (\theta _{i} )} \le \frac{1}{\cos (\varphi )} $ for $i=1,...,n$, the facts that $v_{i} (t)\in (0,v_{\max } )$, $\theta _{i} (t)\in (-\varphi ,\varphi )$ for all $t\ge0$, $i=1,...,n$, and formulas
\[\begin{aligned}
&\dot{F}_{i} (t)=\frac{\sin (\theta _{i} (t))u_{i} (t)}{\cos ^{2} (\theta _{i} (t))} \left(k_{i} (w(t))(v_{i} (t)\cos (\theta _{i} (t))-v^{*} )\right.\\&\left.+\Lambda _{i} (w(t))\right)  \\
&-\frac{1}{\cos (\theta _{i} (t))} \left(\frac{d}{dt} k_{i} (w(t))(v_{i} (t)\cos (\theta _{i} (t))-v^{*} )\right.\\&\left.+\frac{d}{dt} \left(\Lambda _{i} (w(t))\right)\right) \\
&{-\frac{k_{i} (w(t))}{\cos (\theta _{i} (t))} \left(F_{i} (t)\cos (\theta _{i} (t))-v_{i} (t)u_{i} (t)\sin (\theta _{i} (t))\right)} 
\end{aligned}\] 
we conclude that each $\dot{F}_{i} (t)$ is bounded for $i=1,...,n$.

Using \eqref{GrindEQ__2_2_}, \eqref{GrindEQ__3_13_}, \eqref{GrindEQ__5_54_}, \eqref{GrindEQ__5_56_}, \eqref{GrindEQ__5_31_}, \eqref{GrindEQ__5_32_}, \eqref{GrindEQ__5_60_}, and \eqref{GrindEQ__5_61_} we also have that $\frac{d}{dt} \left(Z_{i} (w(t))\right)$ for $i=1,...,n$ are bounded. Combining the fact that $v_{i} (t)\in (0,v_{\max } )$ for $i=1,...,n$ and \eqref{GrindEQ__3_32_}, gives a lower bound for all speeds, i.e., each $\frac{1}{v_{i} (t)} $ is bounded for $i=1,...,n$.

\noindent Let 
\begin{equation} \label{GrindEQ__5_63_} 
\begin{aligned}
Y_{i} (t):=&\left(v^{*} +\frac{A}{v_{i} (t)\left(\cos (\theta _{i} (t))-\cos (\varphi )\right)^{2} } \right.\\
&+v_{i} (t)\cos (\theta _{i} (t))(b-1)\Biggr)^{-1} .   
\end{aligned}
\end{equation} 
Definition \eqref{GrindEQ__5_63_} and \eqref{GrindEQ__5_55_} imply that $Y_{i} (t)$ is bounded for all $i=1,...,n$. Moreover, \eqref{GrindEQ__5_54_}, \eqref{GrindEQ__5_55_}, \eqref{GrindEQ__5_56_}, boundedness of $\frac{1}{v_{i} (t)} $ , inequalities $H(w_{0} )\ge \frac{A}{\cos (\theta _{i} )-\cos (\varphi )} $ for $i=1,...,n$ (recall \eqref{GrindEQ__3_9_} and \eqref{GrindEQ__5_49_}), and formula
\[\begin{aligned} 
&\left(v^{*} +\frac{A}{v_{i} (t)\left(\cos (\theta _{i} (t))-\cos (\varphi )\right)^{2} } +v_{i} (t)\cos (\theta _{i} (t))(b-1)\right)^{2}\\
&\times \frac{d}{d\, t} \left(Y_{i} (t)\right)\\ 
&=\frac{2A\sin (\theta _{i} (t))u_{i} (t)}{\left(\cos (\theta _{i} (t))-\cos (\varphi )\right)^{3} v_{i} (t)} +\frac{AF_{i} (t)}{\left(\cos (\theta _{i} (t))-\cos (\varphi )\right)^{2} v_{i}^{2} (t)} \\ 
&+\left(b-1\right)\sin (\theta _{i} )v_{i} (t)u_{i} (t)-\left(b-1\right)\cos (\theta _{i} (t))F_{i} (t) \end{aligned}\]

\noindent imply that each $\frac{d}{d\, t} \left(Y_{i} (t)\right)$ is bounded for $i=1,...,n$. Finally, since \eqref{GrindEQ__5_50_} holds, it follows that $U'(y_{i} (t))$ and $U''(y_{i} (t))$ are bounded for all $i=1,...,n$ which in conjunction with \eqref{GrindEQ__5_54_}, \eqref{GrindEQ__5_56_}, $v_{i} (t)\in (0,v_{\max } )$, $\theta _{i} (t)\in \left(-\varphi ,\varphi \right)$ for all $t\ge 0$, imply that each $\frac{d}{dt} \left(U'(y_{i} (t))\right)$ is bounded for $i=1,...,n$.

Using \eqref{GrindEQ__5_55_} and the facts that for all $i=1,...,n$, $Y_i(t)$, $\frac{d}{dt} \left(Y_{i} (t)\right)$,  $\frac{d}{d\, t} \left(\sum _{j\ne i}p_{i,j} V'_{i,j} (d_{i,j} (t))\frac{(y_{i} (t)-y_{j} (t))}{d_{i,j} (t)}  \right)$, $F_{i} (t)$, $\dot{F}_{i} (t)$, $\frac{d}{dt} \left(U'(y_{i} (t))\right)$, $u_{i} (t)$, $\frac{d}{dt} \left(Z_{i} (w(t))\right)$,  are bounded, and formulas
\[\begin{aligned}
&\dot{u}_{i} (t)=\\
&\frac{d}{dt} Y_{i} (t)\Biggl(Z_{i} (w(t))-U'(y_{i} (t)) -b\sin (\theta _{i} (t))F_{i} (t) \\
&\left.-\sum _{j\ne i}p_{i,j} V'_{i,j} (d_{i,j} (t))\frac{(y_{i} (t)-y_{j} (t))}{d_{i,j} (t)} \right)\\ 
&+Y_{i} (t)\left(\frac{d}{dt} Z_{i} (w(t))-\frac{d}{dt} U'(y_{i} (t))\right.\\
&\left.-\frac{d}{dt} \left(\sum _{j\ne i}p_{i,j} V'_{i,j} (d_{i,j} (t))\frac{(y_{i} (t)-y_{j} (t))}{d_{i,j} (t)}  \right)\right) \\
& -Y_{i} (t)b\left(\cos (\theta _{i} (t))u_{i} (t)F_{i} (t)+\sin (\theta _{i} (t))\dot{F}_{i} (t)\right)
\end{aligned}\]

\noindent for  we conclude that $\dot{u}_{i} (t)$ are bounded for all $i=1,...,n$. Then, \eqref{GrindEQ__3_33_}, is a direct consequence of Lemma 3. This completes the proof. $\triangleleft $

\section{Conclusions}

In the present work, we have applied a CLF methodology to design nonlinear cruise controllers for the two-dimensional movement of autonomous vehicles on lane-free roads. The CLF were based on measures of the total energy of the system. By expressing the kinetic energy as in Newtonian or relativistic mechanics, two families of controllers were obtained, both guaranteeing that: (i) the vehicles do not collide with each other or with the boundary of the road; (ii) the speeds of all vehicles are always positive and remain below a given speed limit; (iii) all vehicle speeds converge to a given longitudinal speed set-point; and, (iv) the accelerations, lateral speeds, and orientations of all vehicles tend to zero. The proposed families of cruise controllers are decentralized (per vehicle) and require either the measurement only of the distances from adjacent vehicles (inviscid cruise controllers) or the measurement of speeds and distances from adjacent vehicles (viscous cruise controllers). 
 
Finally, we have formally derived the corresponding macroscopic models consisting of a conservation equation and a momentum equation with pressure and viscous terms. We have shown that, by selecting appropriately the parameters of the cruise controllers, we can directly influence the physical characteristics of the ``traffic fluid'', thus creating an artificial fluid that approximates the traffic flow.  In future work, we will study the expected level of approximation of the emerging traffic due to Assumption H as well as how the effects of lateral movement are reflected on the macroscopic model.


\section*{Appendix: Formal Derivation of Macroscopic Models }
\setcounter{equation}{0}
\renewcommand\theequation{A.\arabic{equation}}

Consider the movement of $n$ identical vehicles with total mass $m>0$ on a straight road under the PRCC \eqref{GrindEQ__3_21_}, \eqref{GrindEQ__3_22_} and under the NCC \eqref{GrindEQ__3_11_}, \eqref{GrindEQ__3_12_}  when Assumption H holds. 

Each vehicle has mass $m/n$ and we define the inter-vehicle distance by
\begin{equation}\label{A.1}
s_{i} =x_{i-1} -x_{i} \quad ,\quad i=2,...,n.
\end{equation}
The microscopic model \eqref{GrindEQ__2_2_} under the PRCC \eqref{GrindEQ__3_21_}, \eqref{GrindEQ__3_22_} is given by the following ODEs:
\begin{align}
\dot{x}_{i} =&v_{i} \quad ,\quad i=1,2,...,n\label{A.2}\\
q(v_{1} )\dot{v}_{1} =&-f\left(v_{1} -v^{*} \right)-n\Phi '(ns_{2} )\nonumber\\&+n^{2} K(ns_{2} )(g(v_{2} )-g(v_{1} ))\label{A.3}\\
q(v_{i} )\dot{v}_{i} =&-f\left(v_{i} -v^{*} \right)+n\Phi '(ns_{i} )-n\Phi '(ns_{i+1} )\nonumber\\&+n^{2} K(ns_{i} )(g(v_{i-1} )-g(v_{i} ))\nonumber\\&+n^{2} K(ns_{i+1} )(g(v_{i+1} )-g(v_{i} )),\label{A.4}\\
&\qquad\qquad\textrm{ for } i=2,...,n-1,\nonumber   \end{align}\begin{align}
q(v_{n} )\dot{v}_{n} =&-f\left(v_{n} -v^{*} \right)+n\Phi '(ns_{n} )\nonumber\\&+n^{2} K(ns_{n} )(g(v_{n-1} )-g(v_{n} ))\label{A.5}
\end{align}
\noindent where $f,g\in C^{1} $ with $f(0)=0$, $x\, f(x)>0$ for $x\ne 0$ and $g'(v)>0$ for $v\in \mathbb{R}$ and $q:(0,v_{\max } )\to (0,+\infty )$ given by \eqref{GrindEQ__4_5_} with $v^{*} \in (0,v_{\max } )$ being the longitudinal speed set-point, and $v_{\max } $ the speed limit of the road ($f,g$ are the functions $f_{1} ,g_{1} $ appearing in \eqref{GrindEQ__3_24_}). The state space of system \eqref{A.2}, \eqref{A.3}, \eqref{A.4}, \eqref{A.5} is the set $\Omega =\{\, (x_{1} ,...,x_{n} ,v_{1} ,...,v_{n} )\in \mathbb{R} ^{n} \times \left(0,v_{\max } \right)^{n} :\, n(x_{i} -x_{i+1} )>L\, ,\, i=1,...,n-1\, \}$.  

Consider solutions $(x(t),v(t))\in \mathbb{R} ^{n} \times \mathbb{R} ^{n} $ of the microscopic model \eqref{A.2}, \eqref{A.3}, \eqref{A.4}, \eqref{A.5} with $x(t)=(x_{1} (t),...,x_{n} (t))\in \mathbb{R} ^{n} $ and $v(t)=(v_{1} (t),...,v_{n} (t))\in \mathbb{R} ^{n} $. We assume that for each $t>0$ there exists an interval $I(t)\subseteq \mathbb{R} $ with $\mathop{\lim }\limits_{n\to +\infty } \left(x_{1} (t)\right)=\sup \left(I(t)\right)$, $\mathop{\lim }\limits_{n\to +\infty } \left(x_{n} (t)\right)=\inf \left(I(t)\right)$ and the following property: for each $x\in I(t)$ there exists a sequence of indices $\{\, i_{n} \in \{ 1,...,n\} \, :\, n=3,4,...\, \}$ with $\mathop{\lim }\limits_{n\to +\infty } \left(x_{i_{n} } (t)\right)=x$. Moreover, consider $C^{2} $ density and speed functions $\rho :\bar{\Omega }\to \left(0,\rho _{\max } \right)$, $v:\bar{\Omega }\to \left(0,v_{\max } \right)$, where $\rho _{\max } :=m/L$ and $\bar{\Omega }=\bigcup _{t>0}\{ t\} \times I(t) $, which satisfy the equations
\begin{align}
\rho (t,x_{i} (t))=&\frac{m}{ns_{i} (t)}, \textrm{ for }t>0,\,\, i=2,...,n\label{A.6}\\
v(t,x_{i} (t))=&v_{i} (t) ,\textrm{ for }t>0, \,\,i=1,...,n.\label{A.7}
\end{align}

Notice that by virtue of definition \eqref{A.6} in conjunction with the fact that \eqref{A.1} (which implies that $s_{i} (t)>L/n$ for $i=2,...,n$), we get that $0<\rho (t,x_{i} (t))<m/L=\rho _{\max } $, for $t>0$, $i=2,...,n$.

Using definition \eqref{A.6} and \eqref{A.1}, \eqref{A.2}, we obtain 
\begin{equation}\label{A.8}
\begin{aligned}
\frac{d}{dt} \rho (t,x_{i} (t))=-\rho (t,x_{i} (t))\frac{v_{i-1} (t)-v_{i} (t)}{x_{i-1} (t)-x_{i} (t)},&\\ \textrm{ for } t>0,\,\, i=2,...,n.                 &
\end{aligned}
\end{equation}

\noindent Using the chain rule, we also have from \eqref{A.7} and \eqref{A.2} that 
\begin{equation}\label{A.9}
\begin{aligned}
\frac{d}{dt} \rho (t,x_{i} (t))=\frac{\partial \, \rho }{\partial \, t} (t,x_{i} (t))+\frac{\partial \, \rho }{\partial \, x} (t,x_{i} (t))v(t,x_{i} (t)),&\\
 \textrm{ for }t>0, \,\,i=2,...,n.&
\end{aligned}
\end{equation}
\noindent Thus, we get from \eqref{A.1}, \eqref{A.8} and \eqref{A.9} for all $t>0$, $i=2,...,n$:
\begin{equation}\label{A.10}
\begin{aligned}
&{\frac{\partial \, \rho }{\partial \, t} (t,x_{i} (t))+\frac{\partial \, \rho }{\partial \, x} (t,x_{i} (t))v(t,x_{i} (t))} \\ 
&{=-\rho (t,x_{i} (t))\frac{v(t,x_{i-1} (t))-v(t,x_{i} (t))}{x_{i-1} (t)-x_{i} (t)} } \\ 
&=-\rho (t,x_{i} (t))\left(\frac{\partial v}{\partial x} (t,x_{i} (t))\right.\\&\left.+\frac{1}{s_{i} (t)} \int _{x_{i} (t)}^{x_{i} (t)+s_{i} (t)}\int _{x_{i} (t)}^{l}\frac{\partial ^{2} v}{\partial \, x^{2} } (t,r)drdl  \right)
\end{aligned}
\end{equation}

\noindent Furthermore, assuming that the quantities $\mathop{\max }\limits_{i=2,...,n} \left(ns_{i} (t)\right)$ and $\mathop{\sup }\limits_{t\ge 0,x\in \mathbb{R} } \left(\left|\frac{\partial ^{2} v}{\partial \, x^{2} } (t,x)\right|\right)$ are bounded for $n\ge 2$ and for all $t>0$, we obtain from \eqref{A.10} as $n\to +\infty $ the continuity equation \eqref{GrindEQ__4_3_}. 

Next, using \eqref{A.7} and \eqref{A.4}, we obtain
\begin{equation*}
\begin{aligned}
&q(v(t,x_{i} (t)))\frac{d}{dt} v(t,x_{i} (t))=\\&-f\left(v(t,x_{i} (t))-v^{*} \right)+n\Phi '(ns_{i} (t))-n\Phi '(ns_{i+1} (t)) \\ 
&+n^{2} K(ns_{i} (t))(g(v_{i-1} (t))-g(v_{i} (t)))\\&+n^{2} K(ns_{i+1} (t))(g(v_{i+1} (t))-g(v_{i} (t)))\\
&\qquad\qquad\textrm{ for }t>0,\,\,i=2,...,n-1.
\end{aligned}
\end{equation*}

\noindent We also obtain from the chain rule
\begin{equation*}
\frac{d}{dt} v(t,x_{i} (t))=\frac{\partial \, v}{\partial \, t} (t,x_{i} (t))+\frac{\partial \, v}{\partial \, x} (t,x_{i} (t))v(t,x_{i} (t)),
\end{equation*}
\noindent The equations above imply that
\begin{equation}\label{A.11}
\begin{aligned}
q(v(t,&x_{i} (t)))\frac{\partial \, v}{\partial \, t} (t,x_{i} (t))\\+q(&v(t,x_{i} (t)))v(t,x_{i} (t))\frac{\partial \, v}{\partial \, x} (t,x_{i} (t))= \\ 
&-f\left(v(t,x_{i} (t))-v^{*} \right)+n\Phi '(ns_{i} (t))-n\Phi '(ns_{i+1} (t)) \\ 
&+n^{2} K(ns_{i} (t))(g(v_{i-1} (t))-g(v_{i} (t)))\\&+n^{2} K(ns_{i+1} (t))(g(v_{i+1} (t))-g(v_{i} (t)))\\
&\qquad\qquad \textrm{ for } t\ge0,\,\, i=2,\ldots,n-1.
\end{aligned}
\end{equation} 
\noindent  Notice that definition \eqref{A.6} implies that 
\begin{equation}\label{A.12}
\begin{aligned}
&{s_{i} (t)-s_{i+1} (t)=\frac{m}{n} \left(\frac{1}{\rho (t,x_{i} (t))} -\frac{1}{\rho (t,x_{i+1} (t))} \right)} \\
&=-\frac{m^{2} }{n^{2} \rho (t,x_{i+1} (t))\rho ^{2} (t,x_{i} (t))} \frac{\partial \, \rho }{\partial \, x} (t,x_{i} (t))\\
&+\frac{m}{n} \int _{x_{i} (t)}^{x_{i} (t)-s_{i+1} (t)}\int _{x_{i} (t)}^{l}\frac{1}{\rho ^{2} (t,r)} \frac{\partial ^{2} \, \rho }{\partial \, x^{2} } (t,r)dr dl  \\ 
&{-\frac{2m}{n} \int _{x_{i} (t)}^{x_{i} (t)-s_{i+1} (t)}\int _{x_{i} (t)}^{l}\frac{1}{\rho ^{3} (t,r)} \left(\frac{\partial \, \rho }{\partial \, x} (t,r)\right)^{2} dr dl } 
\end{aligned}
\end{equation}

\noindent Combining \eqref{A.6}, \eqref{A.12}, we get:
\begin{equation}\label{A.13}
\begin{aligned}
n\Phi '&\left(ns_{i} (t)\right)-n\Phi '\left(ns_{i+1} (t)\right) \\ 
=&n^{2} \Phi ''(ns_{i} (t))\left(s_{i} (t)-s_{i+1} (t)\right)\\&-n^{3} \int _{s_{i+1} (t)}^{s_{i} (t)}\int _{l}^{s_{i} (t)}\Phi '''(nr)dr dl 
\end{aligned}
\end{equation}
\noindent Next, define 
\begin{align}
\varphi (s):=&K(s), \textrm{ for } s>L,\label{A.14}\\
\theta (s):=&s^{2} K(s), \textrm{ for } s>L, \label{A.15}\\
w(t,x):=&g(v(t,x)),\textrm{ for }t>0, \,\,x\in \mathbb{R} .\label{A.16}
\end{align}

\noindent By virtue of \eqref{A.6}, \eqref{A.14}, \eqref{A.15}, we get: 
\begin{equation}\label{A.17}
\begin{aligned}
&n^{2} K(ns_{i} (t))\left(w(t,x_{i-1} (t))-w(t,x_{i} (t))\right)\\
&\qquad+n^{2} K(ns_{i+1} (t))\left(w(t,x_{i+1} (t))-w(t,x_{i} (t))\right) \\ 
&=n^{2} K(ns_{i} (t))\left(\frac{\partial \, w}{\partial \, x} (t,x_{i} (t))s_{i} (t)+\frac{1}{2} \frac{\partial ^{2} \, w}{\partial \, x^{2} } (t,x_{i} (t))s_{i}^{2} (t)\right.\\
&\qquad\left.+\int _{x_{i} (t)}^{x_{i} (t)+s_{i} (t)}\int _{x_{i} (t)}^{l}\int _{x_{i} (t)}^{r}\frac{\partial ^{3} \, w}{\partial \, x^{3} } (t,\xi )d\xi  dr dl \right) \\
&+n^{2} K(ns_{i+1} (t))\left(-\frac{\partial \, w}{\partial \, x} (t,x_{i} (t))s_{i+1} (t)\right.\\&\left.+\frac{1}{2} \frac{\partial ^{2} \, w}{\partial \, x^{2} } (t,x_{i} (t))s_{i+1}^{2} (t)\right.\\
&\qquad\left.+\int _{x_{i} (t)}^{x_{i} (t)-s_{i+1} (t)}\int _{x_{i} (t)}^{l}\int _{x_{i} (t)}^{r}\frac{\partial ^{3} \, w}{\partial \, x^{3} } (t,\xi )d\xi  dr dl \right) \\
&=n^{2} \frac{\partial \, w}{\partial \, x} (t,x_{i} (t))\varphi '\left(ns_{i} (t)\right)(s_{i} (t)-s_{i+1} (t))\\
&+\theta (ns_{i} (t))\frac{\partial ^{2} \, w}{\partial \, x^{2} } (t,x_{i} (t)) \\
&+n^{2} K(ns_{i} (t))\int _{x_{i} (t)}^{x_{i} (t)+s_{i} (t)}\int _{x_{i} (t)}^{l}\int _{x_{i} (t)}^{r}\frac{\partial ^{3} \, w}{\partial \, x^{3} } (t,\xi )d\xi  dr dl\\
&  +n^{2} K(ns_{i+1} (t))\int _{x_{i} (t)}^{x_{i} (t)-s_{i+1} (t)}\int _{x_{i} (t)}^{l}\int _{x_{i} (t)}^{r}\frac{\partial ^{3} \, w}{\partial \, x^{3} } (t,\xi )d\xi  dr dl \\ 
&+n\frac{\partial \, w}{\partial \, x} (t,x_{i} (t))\int _{ns_{i+1} (t)}^{ns_{i} (t)}\int _{ns_{i} (t)}^{l}\varphi '(r)dr dl \\
&+\frac{1}{2} \frac{\partial ^{2} \, w}{\partial \, x^{2} } (t,x_{i} (t))\int _{ns_{i} (t)}^{ns_{i+1} (t)}\theta '(l)dl 
\end{aligned}
\end{equation}

Using \eqref{A.6}, \eqref{A.12}, \eqref{A.13}, \eqref{A.17} and assuming that the quantities $\mathop{\max }\limits_{i=2,...,n} \left(ns_{i} (t)\right)$ and $\mathop{\sup }\limits_{t\ge 0,x\in \mathbb{R} } \left|\frac{\partial ^{2} \, \rho }{\partial \, x^{2} } (t,x)\right|$,\linebreak $\mathop{\sup }\limits_{t\ge 0,x\in \mathbb{R} } \left|\frac{\partial \, \rho }{\partial \, x} (t,x)\right|$, $\mathop{\sup }\limits_{t\ge 0,x\in \mathbb{R} } \frac{1}{\rho (t,x)} $, $\mathop{\sup }\limits_{t\ge 0,x\in \mathbb{R} } \left|\Phi ''\frac{m}{\rho (t,x)} \right|$,\linebreak $\mathop{\sup }\limits_{t\ge 0,x\in \mathbb{R} } \left|\Phi '''\frac{m}{\rho (t,x)} \right|$, $\mathop{\sup }\limits_{t\ge 0,x\in \mathbb{R} } \left|\frac{\partial \, w}{\partial \, x} (t,x)\right|$, $\mathop{\sup }\limits_{t\ge 0,x\in \mathbb{R} } \left|\frac{\partial ^{2} \, w}{\partial \, x^{2} } (t,x)\right|$,\linebreak $\mathop{\sup }\limits_{t\ge 0,x\in \mathbb{R} } \left|\frac{\partial ^{3} \, w}{\partial \, x^{3} } (t,x)\right|$, $\mathop{\sup }\limits_{t\ge 0,x\in \mathbb{R} } \left|\varphi '\frac{m}{\rho (t,x)} \right|$, $\mathop{\sup }\limits_{t\ge 0,x\in \mathbb{R} } \left|\theta '\frac{m}{\rho (t,x)} \right|$, \linebreak$\mathop{\sup }\limits_{t\ge 0,x\in \mathbb{R} } K\frac{m}{\rho (t,x)} $ are bounded for $n\ge 2$ and for all $t>0$, we get:
\begin{align}
&s_{i} (t)-s_{i+1} (t)=\nonumber\\&-\frac{m^{2} }{n^{2} \rho (t,x_{i+1} (t))\rho ^{2} (t,x_{i} (t))} \frac{\partial \, \rho }{\partial \, x} (t,x_{i} (t))+O(n^{-3} )=O(n^{-2})\label{A.18}\\
&{n\Phi '\left(ns_{i} (t)\right)-n\Phi '\left(ns_{i+1} (t)\right)} \nonumber\\ 
&=-\Phi ''\left(\frac{m}{\rho (t,x_{i} (t))} \right)\frac{m^{2} }{\rho (t,x_{i+1} (t))\rho ^{2} (t,x_{i} (t))} \frac{\partial \, \rho }{\partial \, x} (t,x_{i} (t))\nonumber\\&+O(n^{-1} )  \label{A.19}\end{align}\begin{align}
&\begin{aligned}
&n^{2} K(ns_{i} (t))\left(w(t,x_{i-1} (t))-w(t,x_{i} (t))\right)\\&+n^{2} K(ns_{i+1} (t))\left(w(t,x_{i+1} (t))-w(t,x_{i} (t))\right) \\ 
&=-\varphi '\left(\frac{m}{\rho (t,x_{i} (t))} \right)\frac{m^{2} }{\rho (t,x_{i+1} (t))\rho ^{2} (t,x_{i} (t))} \\&\times\frac{\partial \, w}{\partial \, x} (t,x_{i} (t))\frac{\partial \, \rho }{\partial \, x} (t,x_{i} (t)) \\ 
&{+\theta \left(\frac{m}{\rho (t,x_{i} (t))} \right)\frac{\partial ^{2} \, w}{\partial \, x^{2} } (t,x_{i} (t))+O\left(n^{-1} \right)}
\end{aligned}\label{A.20}
\end{align} 

\noindent Therefore, we obtain the following equation we obtain from \eqref{A.11}, \eqref{A.19}, \eqref{A.20} as $n\to +\infty $ :
\begin{equation}\label{A.21}
\begin{aligned}
q(v)\frac{\partial \, v}{\partial \, t} +q(v)v\frac{\partial \, v}{\partial \, x}& =-f\left(v-v^{*} \right)-\Phi ''\left(\frac{m}{\rho } \right)\frac{m^{2} }{\rho ^{3} } \frac{\partial \, \rho }{\partial \, x}\\& -\varphi '\left(\frac{m}{\rho } \right)\frac{m^{2} }{\rho ^{3} } \frac{\partial \, w}{\partial \, x} \frac{\partial \, \rho }{\partial \, x} +\theta \left(\frac{m}{\rho } \right)\frac{\partial ^{2} \, w}{\partial \, x^{2} }
\end{aligned}
\end{equation}

\noindent Defining 
\begin{equation}\label{A.22}
P(\rho ):=z-m\Phi '\left(\frac{m}{\rho } \right), \,\,\,\mu (\rho ):=\frac{m^{2} }{\rho } K\left(\frac{m}{\rho } \right)
\end{equation}
\noindent where $z\in \mathbb{R} $ is an arbitrary constant and using definitions \eqref{A.14}, \eqref{A.15}, \eqref{A.16}, we obtain from \eqref{A.21} equation \eqref{GrindEQ__4_4_}. Conditions \eqref{GrindEQ__4_10_}, \eqref{GrindEQ__4_1_} are direct consequences of definitions \eqref{A.22}, the facts that $\mathop{\lim }\limits_{d\to L^{+} } \left(\Phi (d)\right)=+\infty $, $\Phi (d)=K(d)=0$ for $d\ge \lambda $ and definitions $\rho _{\max } :=m/L$, $\bar{\rho }:=m/\lambda $.

The microscopic model \eqref{GrindEQ__2_2_} under the NCC \eqref{GrindEQ__3_11_}, \eqref{GrindEQ__3_12_} is given by the ODEs \eqref{A.2} and the following ODEs:
\begin{align}
\dot{v}_{1} =&-(\gamma +h(G_{1} )\left(v_{1} -v^{*} \right)+G_{1} \label{A.23}\\
\dot{v}_{i} =&-(\gamma +h(G_{i} )\left(v_{i} -v^{*} \right)+G_{i} ,  \textrm{ for }i=2,...,n-1\label{A.24}\\
\dot{v}_{n} =&-(\gamma +h(G_{n} ))\left(v_{n} -v^{*} \right)+G_{n} \label{A.25}
\end{align}
\noindent where 
\begin{equation}\label{A.26}
h(s)=\frac{v_{\max } r(s)}{v^{*} (v_{\max } -v^{*} )} -\frac{s}{v^{*} }
\end{equation} 
\noindent and
\begin{equation}\label{A.27}
\begin{aligned}
{G_{1} } & {=}  {-n\Phi '(ns_{2} )+n^{2} K(ns_{2} )(g(v_{2} )-g(v_{1} ))} \\
 {G_{i} } & {=}  {n\Phi '(ns_{i} )-n\Phi '(ns_{i+1} )+n^{2} K(ns_{i} )(g(v_{i-1} )-g(v_{i} ))} \\ 
& { +n^{2} K(ns_{i+1} )(g(v_{i+1} )-g(v_{i} )),\, \, \, \, \, \, i}  {=}  {1,...,n} \\
  {G_{n} } & {=}  {n\Phi '(ns_{n} )+n^{2} K(ns_{n} )(g(v_{n-1} )-g(v_{n} ))} \end{aligned}
\end{equation} 

\noindent with $\gamma >0$, $g'(v)>0$ for $v\in \mathbb{R}$, ($g$ is the function $g_{1} $ appearing in \eqref{GrindEQ__3_24_}), and $r$ satisfying \eqref{GrindEQ__3_16_}. The state space is $\Omega =\{\, (x_{1} ,...,x_{n} ,v_{1} ,...,v_{n} )\in \mathbb{R} ^{n} \times \left(0,v_{\max } \right)^{n} :\, n(x_{i} -x_{i+1} )>L\, ,\, i=1,...,n-1\, \}$.  

Consider solutions $(x(t),v(t))\in \mathbb{R} ^{n} \times \mathbb{R} ^{n} $ of the microscopic model \eqref{A.2}, \eqref{A.23}, \eqref{A.24}, \eqref{A.25} with $x(t)=(x_{1}(t),...$  $,x_{n} (t))\in \mathbb{R}^{n}$ and $v(t)=(v_{1} (t),\ldots,v_{n} (t))\in \mathbb{R} ^{n} $. We assume that for each $t>0$ there exists an interval $I(t)\subseteq \mathbb{R} $ with $\mathop{\lim }\limits_{n\to +\infty } \left(x_{1} (t)\right)=\sup \left(I(t)\right)$, $\mathop{\lim }\limits_{n\to +\infty } \left(x_{n} (t)\right)=\inf \left(I(t)\right)$ and the following property: for each $x\in I(t)$ there exists a sequence of indices $\{\, i_{n} \in \{ 1,...,n\} \, :\, n=3,4,...\, \}$ with $\mathop{\lim }\limits_{n\to +\infty } \left(x_{i_{n} } (t)\right)=x$. Moreover, consider $C^{2} $ density and speed functions $\rho :\bar{\Omega }\to \left(0,\rho _{\max } \right)$, $v:\bar{\Omega }\to \left(0,v_{\max } \right)$, a $C^{0} $ function $G:\bar{\Omega }\to \left(0,\rho _{\max } \right)$, where $\rho _{\max } :=m/L$ and $\bar{\Omega }=\bigcup _{t>0}\{ t\} \times I(t) $, which satisfy equations \eqref{A.6}, \eqref{A.7} and 
\begin{equation}\label{A.28}
G(t,x_{i} (t))=G_{i} (t), \textrm{ for } t\ge0,\,\,i=1,...,n
\end{equation}

\noindent Using \eqref{A.7}, \eqref{A.24}, \eqref{A.25}, \eqref{A.26}, and the chain rule we obtain
\begin{equation}\label{A.29}
\begin{aligned}
\frac{\partial \, v}{\partial \, t} (t,x_{i} (t))&+v(t,x_{i} (t))\frac{\partial \, v}{\partial \, x} (t,x_{i} )=G_{i} (t), \\
&-(\gamma +h(G_{i} (t)))(v(t,x_{i} (t))-v^{*} )\\
&\textrm{ for } t>0,\,\, i=2,...,n-1.  
\end{aligned}
\end{equation}

\noindent Then, by \eqref{A.12}, \eqref{A.13}, \eqref{A.17}, and definitions \eqref{A.14}, \eqref{A.15}, \eqref{A.16}, and \eqref{A.27}, we get
\begin{equation}\label{A.30}
\begin{aligned}
&{G_{i} (t)=n\Phi '(ns_{i} )-n\Phi '(ns_{i+1} )} \\ 
&+n^{2} K(ns_{i} (t))\left(w(t,x_{i-1} (t))-w(t,x_{i} (t))\right)\\
&+n^{2} K(ns_{i+1} (t))\left(w(t,x_{i+1} (t))-w(t,x_{i} (t))\right) \\ 
&=n^{2} \Phi ''(ns_{i} (t))\left(s_{i} (t)-s_{i+1} (t)\right)\\
&-n^{3} \int _{s_{i+1} (t)}^{s_{i} (t)}\int _{l}^{s_{i} (t)}\Phi '''(nr)dr dl  \\
&+n^{2} \frac{\partial \, w}{\partial \, x} (t,x_{i} (t))\varphi '\left(ns_{i} (t)\right)(s_{i} (t)-s_{i+1} (t))\\
&+\theta (ns_{i} (t))\frac{\partial ^{2} \, w}{\partial \, x^{2} } (t,x_{i} (t)) \\
&+n^{2} K(ns_{i} (t))\int _{x_{i} (t)}^{x_{i} (t)+s_{i} (t)}\int _{x_{i} (t)}^{l}\int _{x_{i} (t)}^{r}\frac{\partial ^{3} \, w}{\partial \, x^{3} } (t,\xi )d\xi  dr dl\\
& +n^{2} K(ns_{i+1} (t))\int _{x_{i} (t)}^{x_{i} (t)-s_{i+1} (t)}\int _{x_{i} (t)}^{l}\int _{x_{i} (t)}^{r}\frac{\partial ^{3} \, w}{\partial \, x^{3} } (t,\xi )d\xi  dr dl \\
&+n\frac{\partial \, w}{\partial \, x} (t,x_{i} (t))\int _{ns_{i+1} (t)}^{ns_{i} (t)}\int _{ns_{i} (t)}^{l}\varphi '(r)dr dl \\
&+\frac{1}{2} \frac{\partial ^{2} \, w}{\partial \, x^{2} } (t,x_{i} (t))\int _{ns_{i} (t)}^{ns_{i+1} (t)}\theta '(l)dl 
\end{aligned}
\end{equation}

\noindent Using \eqref{A.6}, \eqref{A.12}, \eqref{A.30}, and assuming that the quantities $\mathop{\max }\limits_{i=2,...,n} \left(ns_{i} (t)\right)$ and $\mathop{\sup }\limits_{t\ge 0,x\in \mathbb{R} } \left|\frac{\partial ^{2} \, \rho }{\partial \, x^{2} } (t,x)\right|$, $\mathop{\sup }\limits_{t\ge 0,x\in \mathbb{R} } \left|\frac{\partial \, \rho }{\partial \, x} (t,x)\right|$, $\mathop{\sup }\limits_{t\ge 0,x\in \mathbb{R} } \left(\frac{1}{\rho (t,x)} \right)$, $\mathop{\sup }\limits_{t\ge 0,x\in \mathbb{R} } \left|\Phi ''\left(\frac{m}{\rho (t,x)} \right)\right|$, $\mathop{\sup }\limits_{t\ge 0,x\in \mathbb{R} } \left|\Phi '''\left(\frac{m}{\rho (t,x)} \right)\right|$, $\mathop{\sup }\limits_{t\ge 0,x\in \mathbb{R} } \left|\frac{\partial \, w}{\partial \, x} (t,x)\right|$, $\mathop{\sup }\limits_{t\ge 0,x\in \mathbb{R} } \left|\frac{\partial ^{2} \, w}{\partial \, x^{2} } (t,x)\right|$, $\mathop{\sup }\limits_{t\ge 0,x\in \mathbb{R} } \left|\frac{\partial ^{3} \, w}{\partial \, x^{3} } (t,x)\right|$, $\mathop{\sup }\limits_{t\ge 0,x\in \mathbb{R} } \left|\varphi '\left(\frac{m}{\rho (t,x)} \right)\right|$, $\mathop{\sup }\limits_{t\ge 0,x\in \mathbb{R} } \left|\theta '\left(\frac{m}{\rho (t,x)} \right)\right|$, $\mathop{\sup }\limits_{t\ge 0,x\in \mathbb{R} } K\left(\frac{m}{\rho (t,x)} \right)$  are bounded for $n\ge 2$ and for all $t>0$, we get:
\begin{equation}\label{A.31}
\begin{aligned}  &{G_{i} (t)}   {=}  \\& {-\Phi ''\left(\frac{m}{\rho (t,x_{i} (t))} \right)\frac{m^{2} }{\rho (t,x_{i+1} (t))\rho ^{2} (t,x_{i} (t))} \frac{\partial \, \rho }{\partial \, x} (t,x_{i} (t))} \\&-\varphi '\left(\frac{m}{\rho (t,x_{i} (t))} \right)\frac{m^{2} }{\rho (t,x_{i+1} (t))\rho ^{2} (t,x_{i} (t))}\\&\times\frac{\partial \, w}{\partial \, x} (t,x_{i} (t))\frac{\partial \, \rho }{\partial \, x} (t,x_{i} (t)) \\ &{+\theta \left(\frac{m}{\rho (t,x_{i} (t))} \right)\frac{\partial ^{2} \, w}{\partial \, x^{2} } (t,x_{i} (t))+O\left(n^{-1} \right)} \end{aligned} 
\end{equation}
\noindent Taking into account \eqref{A.14}, \eqref{A.15}, \eqref{A.16}, \eqref{A.22}, we obtain from \eqref{A.26}, \eqref{A.28}, \eqref{A.29}, \eqref{A.31}, as $n\to +\infty $, equation \eqref{GrindEQ__4_6_}.

\end{document}